\documentclass[journal]{IEEEtran}

\usepackage{graphics} 
\usepackage{epsfig} 
\usepackage{mathptmx} 
\usepackage{times} 
\usepackage{amsmath,amssymb} 
\usepackage{amssymb}  
\usepackage{graphicx}
\usepackage{todonotes}
\usepackage{tikz}
\usepackage{verbatim}

\usetikzlibrary{decorations.markings}
\usetikzlibrary{arrows}

\newtheorem{theorem}{Theorem}
\newtheorem{example}[theorem]{Example}
\newtheorem{assumption}{Assumption}
\newtheorem{definition}[theorem]{Definition}

\newtheorem{lemma}[theorem]{Lemma}
\newtheorem{corollary}[theorem]{Corollary}

\newtheorem{remark}{Remark}

\DeclareMathOperator{\diag}{diag}
\DeclareMathOperator{\divergence}{div}
\DeclareMathOperator{\Image}{Im}
\DeclareMathOperator{\aufspann}{span}

\usepackage{lipsum}
\usepackage{color}
\usepackage{tikz,pgfplots}

\title{Ensemble Observability of Linear Systems}

\author{Shen Zeng, Steffen Waldherr, Christian Ebenbauer, and Frank Allg\"{o}wer
\thanks{}
\thanks{Section II and III extend results that were presented at the 53rd IEEE Conference on Decision and Control, see \cite{Zeng2014inverse}.}

\thanks{Shen Zeng, Christian Ebenbauer and Frank Allg\"{o}wer are with the Institute for Systems Theory and Automatic Control,
University of Stuttgart, 70550 Stuttgart, Germany, and thank the German Research Foundation (DFG) for financial support of the project within the Cluster of Excellence in Simulation Technology (EXC 310/2) at the University of Stuttgart.
}%
\thanks{Steffen Waldherr is with the Institute for Automation Engineering, Otto-von-Guericke University Magdeburg, 39106 Magdeburg, Germany. 
}
\thanks{Correspondence to {\tt shen.zeng@ist.uni-stuttgart.de}}
}

\begin{document}

\maketitle
\thispagestyle{empty}
\pagestyle{empty}

\begin{abstract}
We address the observability problem for ensembles that are described by probability distributions. The problem is to reconstruct a probability distribution of the initial state from the time-evolution of the probability distribution of the output under a classical finite-dimensional linear system.
We present two solutions to this problem, one based on formulating the problem as an inverse problem and the other one based on reconstructing all the moments of the distribution. The first approach leads us to a connection between the reconstruction problem and mathematical tomography problems. In the second approach we use the framework of tensor systems to describe the dynamics of the moments which leads to a more systems theoretic treatment of the reconstruction problem. Furthermore we show that both frameworks are inherently related. The appeal of having two dual view points, the first being more geometric and the second one being more systems theoretic, is illuminated in several examples of theoretical or practical importance.
\end{abstract}

\begin{IEEEkeywords}
Observability, ensemble control, tomography, moment dynamics, polynomial systems
\end{IEEEkeywords}

\section{INTRODUCTION}
\IEEEPARstart{T}{he} classical question of observability asks whether it is possible to reconstruct the initial state of a finite-dimensional system via the knowledge of the evolution of the outputs, as introduced by R.\ E.\ Kalman in  \cite{Kalman1959_general_theory}, see also \cite{Kalman1963_mathematical_linear_systems}. 
 The concept of observability has become one of the fundamental concepts of modern control theory.
In this paper we address a novel yet very natural extension to this problem in which we move from points in finite-dimensional space to probability distributions.
We ask ourselves under which conditions it is possible to reconstruct a non-parametric probability distribution of initial states when given only knowledge of the evolution of the probability distribution of outputs.
This basic and fundamental problem is not only of interest in its own right but can also be considered a theoretical foundation for state estimation problems for so-called ensembles appearing in a variety of different fields.

In many fields such as process engineering, cell biology, or quantum systems, one encounters large populations of systems that are governed by the same dynamical process, but which have quantities that are distributed among the population, and which can only be manipulated or observed as a whole. The states of such ensembles are commonly modelled as a density function on the state space of the individual systems.
An example from process engineering are particle systems, where for example polymerization processes are modelled dynamically with a density function over the particle size that evolves dynamically \cite{WangDoy2004,PalisKie2014}.
In cellular biology, populations of heterogeneous cells are described via a density function over the heterogeneous cellular variables, which change dynamically due to cellular physiology \cite{hasenauer2011analysis}. 
The same model class is obtained not only by considering populations of similar individuals, but also by studying a single system with a probabilistic description of model uncertainty. An early example is from chemical kinetics, where this problem has been termed ``stochastic sensitivity analysis'' \cite{CostanzaSei1981}.

These ensembles are now also beginning to attract attention in the control community, as they define a novel setting for classical systems theoretical problems, see e.g. \cite{brockett2010control,brockett2012notes}. Based on their common problem formulation and premise, these problems have been grouped into what is now called ensemble control. Recent work in ensemble control has focused almost exclusively on the control part \cite{li2011ensemble,li2006control, Helmke201469, schoenlein2013necsys,becker2012approximate}. There the fundamental difficulty stems from the premise that all systems in an ensemble receive exactly the same control signal. Therefore, the control of such ensembles is picturesquely called broadcast control (cf. \cite{azuma2013broadcast, wood2008broadcast}) as it can be thought of as one command being sent to all individual systems at the same time via a broadcast. The ensemble observability problem presented here can be considered the natural counterpart of these problems. Just as it is not possible to manipulate systems in the ensemble individually, one also cannot track or observe systems in the ensemble individually.

\begin{figure}[htpb]
\vspace{-0.6cm}
 \centering
    \includegraphics[width=0.43 \textwidth]{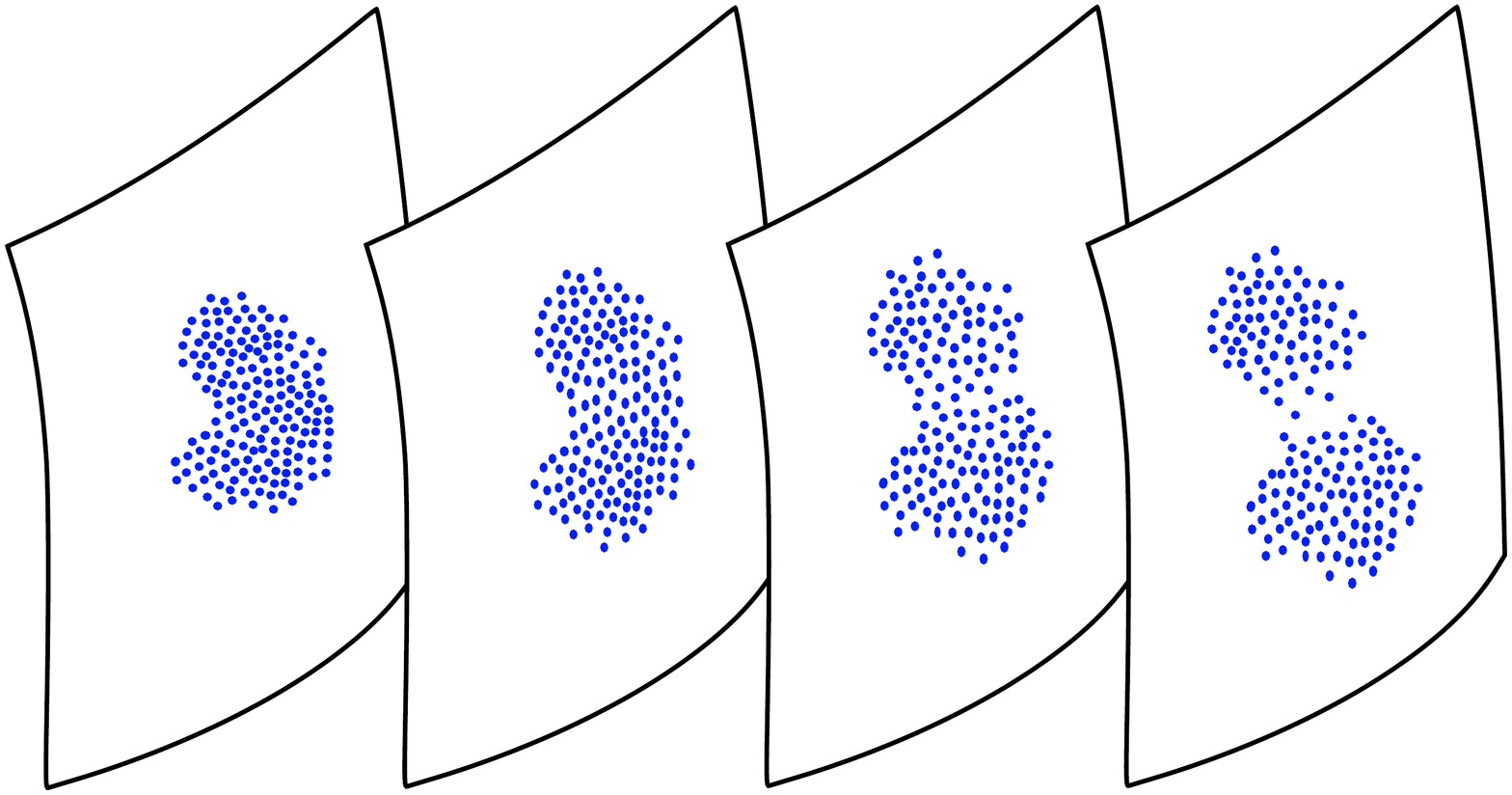} 
	\put(-222,112){\rotatebox{34.3}{Output space}}
	\put(-192,16){\rotatebox{34.3}{$t=t_1$}}
	\put(-142,16){\rotatebox{34.3}{$t=t_2$}}
	\put(-92,16){\rotatebox{34.3}{$t=t_3$}}
	\put(-40,16){\rotatebox{34.3}{$t=t_4$}}
	\vspace{-0.55cm}
      \caption{An illustration of population snapshots. In each time step $t_1, \dots, t_4$ we have a snapshot of certain output values of a population. The crucial point is that in a snapshot, information relating an output value to the individual producing that output value is completely missing.}
\label{fig:aggregate_data}
\end{figure}

An example for such state estimation problem for ensembles is given by the state estimation of heterogeneous cell populations. Due to the experimental circumstances, measurement data in the context of cell populations consists mostly of population snapshots which are provided by high-throughput devices such as flow cytometers, as illustrated in Figure~\ref{fig:aggregate_data}.

While a population snapshot does contain a vast numbers of output measurements, this is at the drawback of losing any information relating output measurements to the individuals that produced them. Thus, while we get a good idea of how the population evolves as a whole, due to the aggregation of data, we cannot tell which point corresponds to which individual. In particular we do not know how a particular individual evolved over time. In fact, state estimation for heterogeneous cell populations portrays a rather drastic example for state estimation problems for ensembles. This is because measuring e.g. protein concentration within a cell often results in killing the cell, making it impossible to measure that cell again. Measurements given in different snapshots therefore stem from different cells within the population. For these reasons, the only adequate mathematical model is to view an output snapshot as a set of \emph{samples} taken from \emph{the output distribution} which evolves according to the dynamics of the structural system. Idealizing the vast number of samples as distributions, this immediately leads to the ensemble observability problem.

State and parameter estimation for heterogeneous cell populations from population snapshots has been recently considered from an applied point of view \cite{hasenauer2010identification,hasenauer2011analysis,zechner2012moment}.
While a range of numerical methods for solving such inverse problems is known (see \cite{banks2012review} for a recent review), the estimation problem has not been characterized from a systems theoretic point of view.
For example, while solutions produced by ad hoc optimization based reconstruction algorithms should fit the measurement data well by construction, it is not at all clear what the precise relations between the solution of the algorithms and the real solution are \cite{banks2012review}. 
Neither are there results establishing the uniqueness of an estimation result a priori.

Besides the introduction and motivation of the ensemble observability problem, our main contribution of this paper is a complete systems theoretic characterization of ensemble observability for linear systems. We apply a measure theoretic description of dynamical systems acting on probability distributions \cite{Lasota1994} and establish a connection of the ensemble observability problem with a classical problem of mathematical tomography and integral geometry \cite{natterer1986mathematics,markoe2006analytic}. This is the problem of reconstructing an internal density from external Radon projections. In fact, as it turns out, the output distributions are Radon projections of the initial distribution. Combining results from the field of mathematical tomography with observability properties of the classical linear system, we obtain sufficient conditions for ensemble observability of continuous ensembles of LTI systems. The connection to mathematical tomography reveals geometric properties of the ensemble observability problem nicely and furthermore allows us to transfer the well-developed computational reconstruction methods from computed tomography to reconstruction problems for dynamical ensemble systems.

While the relation to tomography does provide us with a full algebraic characterization of ensemble observability, for a given LTI system the condition is hard to check in general. Our second approach, which is based on the reconstruction of the moments of the distribution, resolves this problem. In this approach, the dynamics of the moments of the state distribution are governed by systems of homogeneous forms, called tensor systems. Although these systems describe the dynamics of monomials, the systems themselves are linear. Therefore, another sufficient condition for ensemble observability is the observability of every tensor system. In this systems theoretical approach, we can also easily formulate a feasible sufficient condition for ensemble observability of a specific class of single-output systems. As it turns out, the most general condition for ensemble observability is still very restrictive compared to that of classical observability. Lastly, we establish the equivalence of the tomographic and the moment-based approaches, thereby giving a complete picture of the problem.

The organization of this paper is as follows. In Section~\ref{section:problem_formulation} we introduce the ensemble observability problem mathematically and elaborate on the setup. We formulate the direct and inverse problems on the level of probability distributions and probability density functions. By inspecting these formulations from an inverse problems perspective, we identify them as tomography problems. In Section~\ref{section:solution_via_mathematical_tomography} we examine the theory of mathematical tomography and integral geometry which provides results for the uniqueness of the reconstruction problem. We furthermore demonstrate how the tomographic reconstruction methods can be used for the practical reconstruction of initial state distributions. In Section~\ref{section:moment_approach} we pursue an alternative approach based on reconstructing all the moments of the distribution. To this end, we utilize the framework of tensor systems, which describe the dynamics of monomials of the original state. As the dynamics of the monomial systems are linear, the ensemble observability problem can be described via the observability of these tensor systems. We establish, as one of our main results, the equivalence between the tomographic characterization of ensemble observability and this more systems theoretic characterization. The appeal of having a complete framework originating from two different view points is illustrated in several examples.

\section{Ensemble observability}
\label{section:problem_formulation}
To capture the essence of those inverse problems described in the introduction, we consider an observability problem for a finite-dimensional linear time-invariant system 
\begin{align}
\label{classical_linear_system}
\begin{split}
    \dot{x}(t) &= A x(t) , \;\; x(0)=x_0, \\
	y(t)   &= C x(t),
\end{split}
\end{align}
in which the initial state $x_0$ is a random vector, that is a multivariate random variable, with probability distribution $\mathbb P_0$. 

Before going on to discussing this setup mathematically, we would like to elaborate more on the setup in view of population models and thereby introduce some terminology. In view of population models or ensembles, one can think of the setup as a description of a continuum of individual systems that have the same dynamics and measurement outputs given by \eqref{classical_linear_system}, but different initial states. We call \eqref{classical_linear_system} the \emph{structural system} of the ensemble, and $\mathbb P_0$ the \emph{initial distribution} which accounts for the specific heterogeneity of the initial states within the population. The fact that the initial state is a random vector leads to the fact that the output $y(t)$ at any given time is also a random vector. Let $\mathbb P_{y(t)}$ denote its distribution, which we shall call \emph{output distribution}.

Now, analogously to the classical observability problem, the ensemble observability problem aims at reconstructing the initial distribution $\mathbb P_0$  from the evolution of the output distributions $\mathbb P_{y(t)}$. A linear system, for which it is possible to  \emph{uniquely} reconstruct an initial distribution from the time-evolution of the output distribution will be called ensemble observable.
\begin{definition}[Ensemble observability of linear systems]\\
A linear system is called \emph{ensemble observable} if 
\begin{align*}
        (\forall t \ge 0 \;\; \mathbb P_{y(t)|\mathbb P_0'} = \mathbb P_{y(t)|\mathbb P_0''}) \;\; \Rightarrow \;\; \mathbb P_0' = \mathbb P_0'',
\end{align*}
for all continuous initial distributions $\mathbb P_0'$ and $\mathbb P_0''$.
\label{def:ensemble_observability}
\end{definition}
%

It turns out that aiming for a general result for arbitrary non-parametric probability distributions is impossible. This is why we will need to restrict our attention to more specific classes of distributions for our theoretical results later. 




\subsection{The direct problem in a measure theoretic framework}
Clearly, the study of ensemble observability is an inverse problem, just as the classical controllability and observability problems can be naturally viewed as well, see e.g. \cite{luenberger1969optimization}. Thus, before we proceed with addressing the inverse problem, it is reasonable to first discuss the direct problem, i.e. how the output distribution evolves from the initial distribution under the finite-dimensional LTI system~\eqref{classical_linear_system}. Having established this, we can then investigate to which extent we can use this direct problem to go the opposite direction for the inverse problem. 

Since the output distribution $\mathbb P_{y(t)}$ is by definition the probability distribution of the random vector $y(t)$, and since the output is related to the initial state via
$y(t) = Ce^{At} x_0$, the distribution $\mathbb P_{y(t)}$ is simply recognized as the \emph{push-forward measure} of $\mathbb P_0$ under the mapping $x \mapsto Ce^{At}x$.
\begin{figure}[htp]
\vskip -0.2cm
\centering \includegraphics[width=0.48 \textwidth]{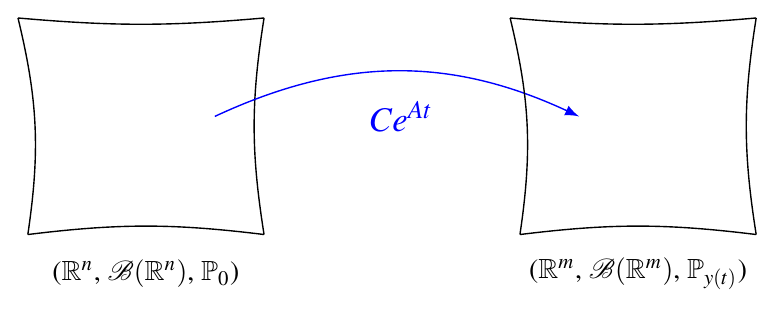} 
\vspace{-0.3cm}
\caption{This figure illustrates the connection between $\mathbb P_0$ and $\mathbb P_{y(t)}$. The distribution $\mathbb P_{y(t)}$ is the push-forward measure of $\mathbb P_0$ under the mapping $Ce^{At}$.}
\label{fig:pushforward}
\end{figure}

This situation is illustrated in Figure \ref{fig:pushforward} and may be formulated as follows. For a measurable set $B_y \subset \mathbb R^m$ we have
\begin{align}
 \mathbb P_{y(t)}(B_y) := \mathbb P_0 ( (Ce^{At})^{-1}(B_y) ),
\label{pushforward}
\end{align}
that is, to compute the probability $\mathbb P_{y(t)}(B_y)$ for a Borel set $B_y \in \mathcal B(\mathbb R^m)$, one pulls back $B_y$ via the mapping $x \mapsto Ce^{At}x$ to obtain the pre-image $(Ce^{At})^{-1}(B_y)$ which is then measured via the probability measure $\mathbb P_0$. 

Assuming furthermore that $\mathbb P_0$ has a probability density function $p_0: \mathbb R^n \to \mathbb R$, we may reformulate \eqref{pushforward} as
\begin{align}
     \mathbb P_{y(t)}(B_y) = \int_{(Ce^{At})^{-1}(B_y)} p_0 \; \text{d}x.
\label{eq:integral_strip}
\end{align}
This establishes the direct problem which can now be used to address the inverse problem of reconstructing $\mathbb P_0$ from $\mathbb P_{y(t)}$.

\subsection{The ensemble observability problem as an inverse problem}
\begin{figure}[htp]
\begin{center}
\vspace{-0.3cm}
 \includegraphics[angle = -90, width=0.34 \textwidth]{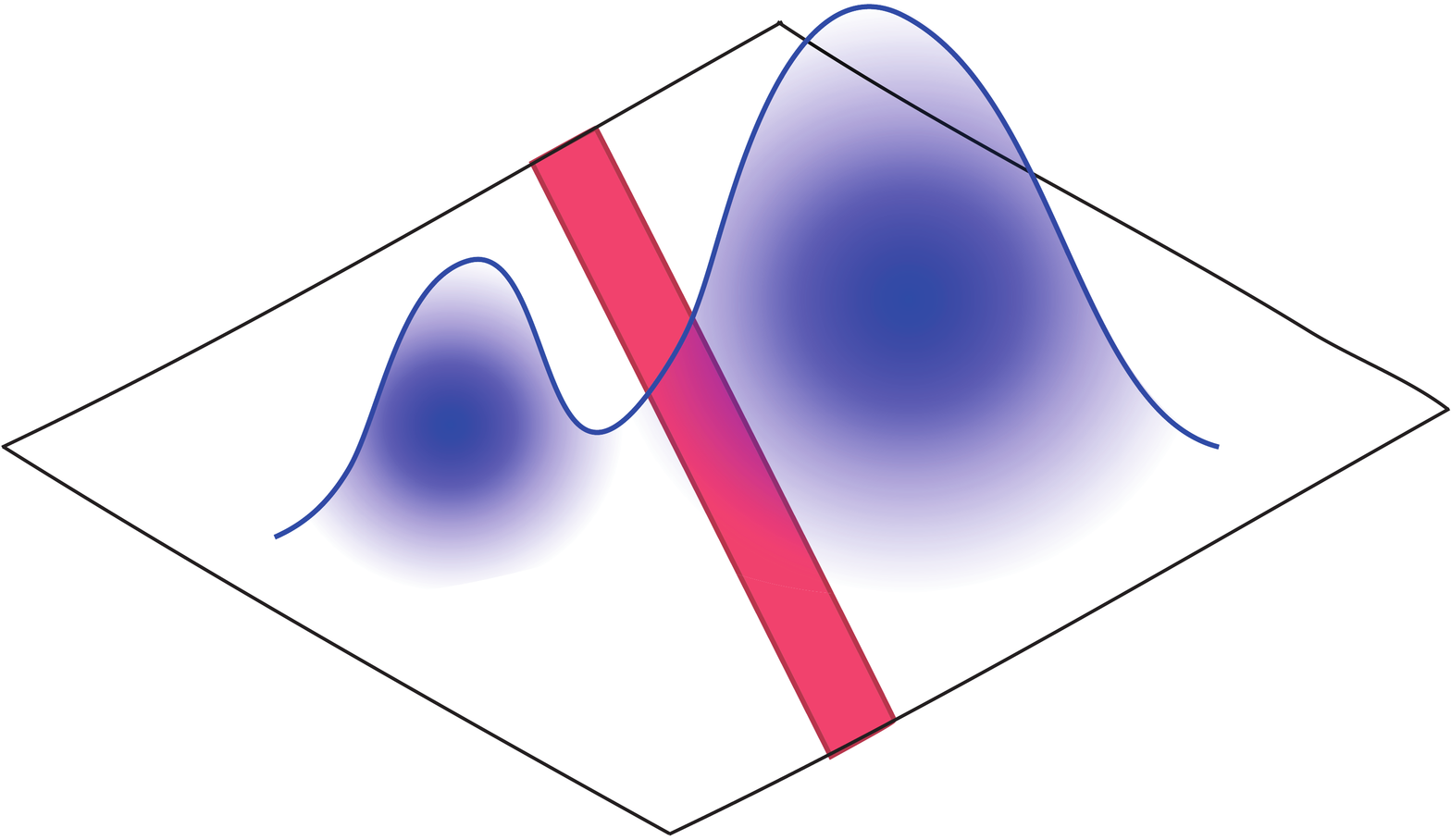} 
  \put (-45,-16) {\textcolor{blue}{unknown $p_0(x)$}}
  \put (-65,-110) {\textcolor{red}{$(Ce^{At})^{-1}(B_y)$}}
\vspace{-0.2cm}
\end{center}
\caption{This figure illustrates the problem that is at the core of our inverse problem: The reconstruction of an unkown density $p_0$ from its integrals along sets $(Ce^{At})^{-1}(B_y)$ for different $t \ge 0$ and $B_y \in \mathcal B(\mathbb R^m)$.}
\label{density_strips}
\end{figure}

Having established the forward relation \eqref{pushforward} between $\mathbb P_0$ and $\mathbb P_{y(t)}$, and the forward relation \eqref{eq:integral_strip} between $p_0$ and $\mathbb P_{y(t)}$, we can already see which information we can infer about the inital distribution from the output distribution for the inverse problem. Since by assumption we know for any time the probability distribution $\mathbb P_{y(t)}$, we now also know the value
\begin{align*}
\mathbb P_0 ( (Ce^{At})^{-1}(B_y) ) = \int_{(Ce^{At})^{-1}(B_y)} p_0 \; \text{d}x
\end{align*}
for all $t \ge 0$ and all $B_y \in \mathcal B(\mathbb R^m)$, which is the probability that an initial state $x_0$ lies in the set $(Ce^{At})^{-1}(B_y)$. We may illustrate the inverse problem as in Figure~\ref{density_strips} now.

Each single piece of information on the initial distribution obtained through the values $\int_{(Ce^{At})^{-1}(B_y)} p_0 \; \text{d}x$ does not provide us with a lot of knowledge about $p_0$ in general, since the output matrix $C \in \mathbb R^{m \times n}$ is in general not injective. Thus in general a set $(Ce^{At})^{-1}(B_y)$ stretches to infinity due to the non-trivial kernel of the matrix $Ce^{At}$. Before studying this inverse problem in full generality, we illustrate in an example, how ensemble observability is related to observability of the structural system mathematically.

\begin{example}
\label{example:illustrative}
We consider the system
\begin{align}
\label{single_cell_example}
 \dot{x}(t) = \begin{pmatrix} -1 & 1 \\ 0 & 0 \end{pmatrix} x(t), \;\; x(0) \sim \mathbb P_0,
\end{align}
see Figure~\ref{fig:flow}, with two different output matrices
\begin{align*}
   C' = \begin{pmatrix} 0 &  1 \end{pmatrix}     \;\; \text{ and } \;\; C'' = \begin{pmatrix} 1 & 0 \end{pmatrix}.
\end{align*}
Concerning the output matrices, we notice that the first output matrix $C'$ leads to an unobservable structural system, while $C''$ renders the structural system observable. 

Returning to our inverse problem, we had already translated the ensemble observability problem to the inverse problem given by (\ref{eq:integral_strip}). It is apparent that the properties of the kernels 
\begin{align*}
\ker Ce^{At}  = e^{-At}(\ker C)
\end{align*} 
will play a crucial role in the reconstruction. In this example we have $\ker C' = \aufspann( \begin{pmatrix} 1 & 0 \end{pmatrix})$ for the first output matrix, which leads to the fact that the sets $(C'e^{At})^{-1}(B_y)$, for a given measurable set $B_y \subset \mathbb R$, are horizontal strips. But by only having at hand integrals over strips that have the same ``orientation'', one cannot uniquely reconstruct $p_0$. 

This is because one can take an arbitrary $p_0$ and shift it along the $x_1$-axis, leaving the resulting output distributions unchanged. More generally, the idea is that due to unobservability, one clearly ends up with a \emph{non-trivial} intersection
\begin{align}
  \bigcap_{t \ge 0} \ker Ce^{At} = \{ x_0 \in \mathbb R^n \; : \; Ce^{At} x_0 \equiv 0 \}.
\label{intersection}
\end{align}
Taking an abritrary density $p_0$ we may shift this density along a \emph{non-zero} vector taken out of the intersection. We end up with two different densities that give the same value when integrated over sets $(Ce^{At})^{-1}(B_y)$. A proof of the general situation will be given in the next subsection.
   \begin{figure}[tbp]
      \centering
       \hspace{0.7cm} \includegraphics[width=0.36 \textwidth]{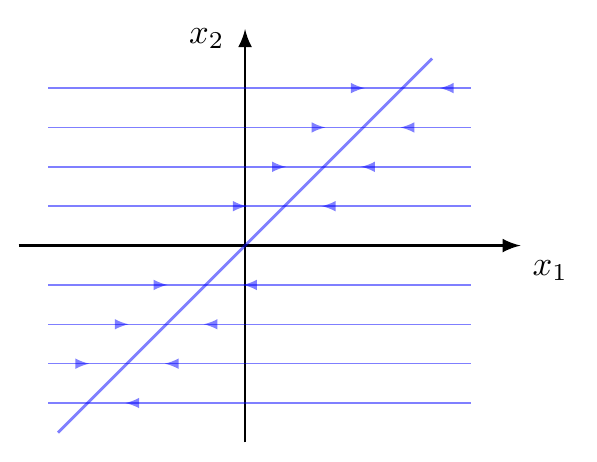} 
\vspace{-0.05cm}
\caption{The phase portrait of the system $\dot{x}(t)=Ax(t)$ given in \eqref{single_cell_example}.}
      \label{fig:flow}
   \end{figure}

It is interesting to see now what will happen in the observable case. There, we find that $\ker C'' = \aufspann( \begin{pmatrix} 0 & 1 \end{pmatrix})$. The evolution of $e^{-At}(\ker C'')$ shows that the kernels are now tilted, which is also to be expected in virtue of a trivial intersection (\ref{intersection}). This is illustrated in Figure \ref{fig:tilted}. Thus, by considering measurements of the output distribution $\mathbb P_{y(t)}$ at different time points, we now obtain integrals of $p_0$ along strips at different angles, thus gaining much more information about the initial density $p_0$ than in the unobservable case. Yet, it is noted that the available range of angles is still constrained by the observability properties of the system. The intriguing question is, whether or not these different pieces of partial information can be put together so as to get full information on $p_0$.

   \begin{figure}[thpb]
 \centering
     \hspace{0.7cm}   \includegraphics[width=0.36 \textwidth]{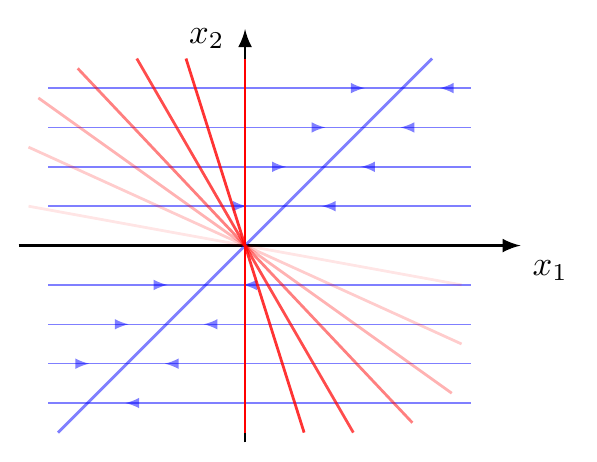} 
\vspace{-0.05cm}
\caption{The evolution of the kernel $\ker C''e^{At}$ is given by transporting the kernel $\ker C''$ with the flow of $\dot{x}(t) = Ax(t)$ back in time.
This results in a tilting of the kernel $\ker C''$, and their intersection is trivial, corresponding to $(A,C'')$ being observable. The color intensities of the kernels indicate their advancement in time.}
\label{fig:tilted}
   \end{figure}
\end{example}

At this point, readers familiar with such inverse problems may have already identified our problem as a \emph{tomography problem}. To further highlight this, we reformulate the inverse problem \eqref{pushforward} just one last time. From the \emph{coarea formula} it follows that the integral in \eqref{eq:integral_strip} is equal to
\begin{align*}
      \int_{B_y} \left( \frac{1}{\sqrt{ \det({Ce^{At}(Ce^{At})^{\top}})}} \int_{(Ce^{At})^{-1}(\{y\})} p_0 \; \text{d}S  \right) \text{d}y,
\end{align*}
provided that the matrix $C$ is not full rank. Here the inner integration with respect to $\text{d}S$ denotes a surface integral over the affine subspaces defined by the equation $Ce^{At}x=y$. It is readily recognized that the density of the output is then
\begin{align}
p_{y(t)}(y) = \frac{1}{\sqrt{ \det({Ce^{At}(Ce^{At})^{\top}})}} \int_{(Ce^{At})^{-1}(\{y\})} p_0 \; \text{d}S.
\label{marginal}
\end{align}
Intuitively, passing from the integral \eqref{density_strips} to the surface integral \eqref{marginal} can be thought of as a concentration of the information $\int_{(Ce^{At})^{-1}(B_y)} p_0 \; \text{d}x$ about $p_0$ by taking the ``width" of a ``strip" as illustrated in Figure \ref{density_strips}, to zero.

To recap, the ensemble observability problem is now recast as the problem of reconstructing a density from its integrals over affine subspaces defined by the equation $Ce^{At}x = y$.
Such problems of reconstructing a density from integrals over affine subspaces, or more generally manifolds, are unique to the theory of mathematical tomography and integral geometry. Our first approach is based on this natural connection. Before we proceed with studying this problem in the tomographic framework, we would like to give several remarks and furthermore give a necessary condition based on the findings of the illustrative example.

\begin{remark}[The case of full state measurements]
\label{full_state}
The formulation \eqref{marginal} is restricted to matrices $C$ that are not full rank. The case in which $C$ is invertible, however, is trivial. There it suffices already to know the output distribution at only one instance $t^{\star}\ge 0$. To see this, assume without loss of generality that $C=I$ such that (\ref{pushforward}) specializes to
\begin{align*}
    \mathbb P_{x(t^\star)}(B) = \mathbb P_0 ( e^{-At^{\star}}(B)).
\end{align*}
To compute for a measureable set $B\in \mathcal B(\mathbb R^n)$ the probability $\mathbb P_0(B)$, we just compute $\mathbb P_{x(t^*)}(e^{At^{\star}} B)$ since
\begin{align*}
 \mathbb P_{x(t^\star)}(e^{At^{\star}} B) = \mathbb P_0 ( e^{-At^{\star}}(e^{At^{\star}} B)) = \mathbb P_0(B)
\end{align*}
due to the invertibility of the matrix exponential. 

Another explanation is via the advection equation or Liouville equation 
\begin{align*}
   \frac{\partial}{\partial t} p(t,x) = - \divergence(p(t,x)F(x)),
\end{align*}
which is a partial differental equation that describes the evolution of a density $x \mapsto p(t,x)$ in a vector field $x \mapsto F(x)$. The important fact to recall here is that the properties of the flow of the vector field are inherited so that the evolution of the density can be described by a flow as well, cf. Frobenius-Perron operators or transfer operators \cite{Lasota1994}. It is noted however, that the measure theoretic explanation is more elementary and that the advection equation in fact originates (mathematically) from the measure theoretic framework, cf. \cite{Lasota1994}.
\end{remark}

\begin{remark}
Having mentioned the Liouville equation, we further note that it is a linear partial differential equation, which suggests that one may study this problem in an infinite-dimensional linear systems framework \cite{curtain1995introduction}. In contrast to the typical output measurements considered in infinite-dimensional linear systems theory, our output is given by
\begin{align*}
  p_{y(t)}(y) = \frac{1}{\sqrt{ \det(CC^{\top})}}  \int_{C^{-1}(\{y\})} p(t,x) \, \text{d}S,
\end{align*}
which is still a surface integral, now taken over the propagated state density. Such integrals are typically not considered in infinite-dimensional linear systems theory, but are in fact a hallmark of tomography problems. Due to this fact, existing results in infinite-dimensional linear systems theory cannot be directly applied to our problem. 
\end{remark}

\begin{remark}[Incorporation of inputs]
Notice that although we formulated the observability problem only for systems \eqref{classical_linear_system} without input, our results do hold for linear systems with input as well. Given 
\begin{align*}
\begin{split}
    \dot{x}(t) &= A x(t) + B u(t) , \;\; x(0)=x_0, \\
	y(t)   &= C x(t),
\end{split}
\end{align*}
with $x_0$ a random vector and a known input function $t \mapsto u(t)$, we can easily deal with the input by considering this in the pre-images of the new mapping $x \mapsto Ce^{At}x + \int_{0}^t  Ce^{A(t-\tau)} B u(\tau) \, \text{d}\tau$, which would not alter the analysis at all. Another way to deal with this is by canceling the convolution integral out by redefining the output. Here it is important that both the input matrix and the input signal be identical for all individuals among the population, so that $\int_{0}^t  Ce^{A(t-\tau)} B u(\tau) \, \text{d}\tau$ is the same for all individuals. This setup thus describes a population that is steered by one single signal in a broadcast manner, which is also the fundamental premise in ensemble control.
\end{remark}

\subsection{Observability of structural system is necessary}
In this subsection, we present a first theoretical result concerning the necessity of observability of the structural system for the corresponding ensemble observability problem to admit a unique solution. This result is a straightforward generalization of our findings in the illustrative example.

\begin{theorem}[Necessary condition \cite{zeng2013identifiability,waldherr2014identifiability}] 
\label{thm:necessity}
Observability of $(A,C)$ is a necessary condition for ensemble observability.
\end{theorem}

\begin{IEEEproof} We show that under the assumption that $(A,C)$ is not observable, there exist initial densities $p_0{\!'} \ne p_0{\!''}$ for which
\begin{align*}
\int_{(Ce^{At})^{-1}(B_y)} p_0' \; \text{d}x = \int_{(Ce^{At})^{-1}(B_y)} p_0'' \; \text{d}x
\end{align*}
for all $t \ge 0$ and $B_y \in \mathcal B(\mathbb R^n)$. In other words, we construct two distinct initial densities $p_0'$ and $p_0''$ which are indistinguishable from the output distributions. To this end, we fix an arbitrary probability density function $p_0{\!'}$. 

Since $(A,C)$ is not observable, the intersection \eqref{intersection}, which is the unobservable subspace, is non-trivial. Therefore we can pick out of this unobservable subspace a non-zero vector $v$, and given that define a second probability density function by
\begin{align*}
   p_0{\!''} (x) := p_0{\!'}(x+v).
\end{align*}
Clearly these two densities are distinct. Furthermore we have
\begin{align*}
\int_{(Ce^{At})^{-1} (B_y )} p_0{\!''}(x) \; \text{d}x  \, &= \, \int_{(Ce^{At})^{-1} (B_y )} p_0{\!'}(x+v) \; \text{d}x  \\
&=  \, \int_{v+(Ce^{At})^{-1} (B_y )} p_0{\!'}(x) \; \text{d}x
\end{align*}
for all $t \ge 0$ and $B_y \in \mathcal B(\mathbb R^m)$. Lastly, we observe that
\begin{align*}
v+(Ce^{At})^{-1} (B_y ) = (Ce^{At})^{-1} (B_y ),
\end{align*}
since $v \in \ker Ce^{At}$ for all $t \ge 0$. This yields the claim.
\end{IEEEproof}
This result generalizes our finding in the illustrative example where an unobservable structural system lead to the existence of indistinguishable initial distributions in the ensemble observability problem. For the observable structural system in the example on the other hand, it seemed plausible that knowing integrals of the unknown density along infinitely many directions could be sufficient for ensemble observability. Yet, the situation is not as clear as for the necessary condition. In the next section, we work towards a resolution to this problem based on the connection to integral geometry and mathematical tomography.

\section{Solution via mathematical tomography}
\label{section:solution_via_mathematical_tomography}

In this section we first give a brief background on the theory of mathematical tomography and afterwards introduce the mathematical framework. Given this framework, we proceed towards a solution of the ensemble observability problem.

\subsection{Background on mathematical tomography}
Classical tomography can be described as a way to determine the internal structure of an object without having to open it up. Probably the best known example for a tomographic problem is computed tomography. Computed tomography is used for providing cross-sections of e.g. a part of a body for medical diagnosis and is based on the physical properties of an \mbox{X-Ray} beam. Let an X-Ray beam $L$ passing through an object with density $f(x)$ be parameterized by $t$, then the intensity $I(t)$ along the X-Ray beam is attenuated according to the Beer-Lambert law \cite{markoe2006analytic}
\begin{align*}
   \frac{d}{dt} I(t) = - f(L(t)) I(t).
\end{align*}
By virtue of this law, measuring the intensity $I_1 = e^{- \int_{L} f(x) \; \text{d}S} I_0$ of the beam after it went through the object and comparing it with the intensity $I_0$ at which it was emitted, we can compute the value 
\begin{align*}
\int_{L} f(x) \; \text{d}S = \log \Big(\frac{I_0}{I_1}\Big).
\end{align*}
A.\ M.\ Cormack, one of the inventors of computed tomography identified this mathematical problem of reconstructing $f$ from its line integrals \cite{cormack1963representation}, \cite{cormack1964representation}. He proposed a reconstruction method for which he was awarded the Nobel prize in medicine and physiology in 1969 jointly with G.\ N.\ Hounsfield. Only later was it discovered that the mathematical problem was solved already 50 years earlier by mathematician J.\ Radon \cite{radon20051}, for purely mathematical reasons.

The general problem of tomography is the reconstruction of a function from its \emph{Radon transform} \cite{markoe2006analytic}, which in its classic form is a transformation that maps a two-dimensional scalar function $f$ to the transform $Rf$ which is defined on lines $L$, i.e. $Rf(L) = \int_{L} f(x) \; \text{d}S$. This is generalized to the $n$-dimensional case as follows, cf. \cite{markoe2006analytic}.
\begin{definition}[Radon transform]
The Radon transform maps an integrable function $f \in L^1(\mathbb R^n, \mathbb R)$ to its transform \mbox{$Rf: \mathbb S^{n-1} \times \mathbb R \to \mathbb R$} which is given as
\begin{align*}
   Rf(\omega,p) = \int_{ \{x \in \mathbb R^n \; : \; \langle \omega, x \rangle = p \} }  f(x) \; \text{d}S
\end{align*}
whenever the integral exists. Furthermore the function $R_\omega f$ given by $(R_\omega f)(p) = Rf(\omega,p)$ is called Radon projection along $\omega^\perp$, or Radon projection orthogonal to $\omega$.
\end{definition}

At this point let us note the immediate connection to our problem \eqref{marginal}. We recognize that the output densities are nothing but some kind of Radon projections of the initial density along $\ker Ce^{At}$. Thus, the inverse problem of reconstructing the initial density from the output densities is clearly a tomography problem.

\begin{remark}
So far we only introduced the classical transform for hyperplanes, while the dimension of $\ker Ce^{At}$ need not be $n-1$. Before presenting results in full generality however, we would like to further illustrate theoretical reconstruction results only for the Radon transform defined for hyperplanes. This is for brevity of presentation as the mathematical framework for transforms for subspaces of arbitrary dimension becomes more involved. Luckily, in our probabilistic framework we may later take a slightly different route for the general problem, leading to a much shorter derivation.
\end{remark}

Let us now turn towards the solution of the reconstruction of a density from its Radon transforms. The key to the theoretical inversion of the Radon transform lies in its connection to the Fourier transform. If we define the $n$-dimensional Fourier transform as
\begin{align*}
  (\mathcal F_n f)(\xi) =  \int_{\mathbb R^n} f(x) e^{-i \langle x, \xi \rangle} \, \text{d}x,
\end{align*}
then, the connection is given as follows.
\begin{theorem}[Projection Slice Theorem, cf. \cite{markoe2006analytic}]
Consider an integrable function $f \in L^1 (\mathbb R^n,\mathbb R)$ and let $\omega$ be a unit vector. Then we have
\begin{align*}
   (\mathcal F_1 R_{\omega}f) (\sigma) = \mathcal F_n f(\sigma \omega).
\end{align*}
\end{theorem}

That is, the one-dimensional Fourier transform of the Radon projection along $\omega^\perp$ is equal to the  $n$-dimensional Fourier transform of the density restricted to the ``slice'' $\sigma \omega$. As for the Fourier transform, we know that it is a bijection. Therefore, if we have for an integrable function, the Radon projections for all ``directions'' $\omega$, then we know the $n$-dimensional Fourier transform of $f$. By bijectivity, we know $f$ and are done. This is the classical solution to the tomography problem. Unfortunately this result does not apply directly to our problem. 
In contrast to the problem in computed tomography, we may not freely choose the directions at which we can gather Radon projections, but the directions $\ker Ce^{At}$ are inherently determined by the observability properties of the finite-dimensional system as was illustrated in Example~\ref{example:illustrative}.

\subsection{Sufficient conditions for ensemble observability}
For the sufficient condition we draw to the mathematical tomography approach introduced before. To deal with the general case of subspaces of arbitrary dimension, and the limited direction problem, we turn towards a probabilistic approach to the tomography problem. 

First of all, we would like to reformulate the Projection Slice Theorem in the probabilistic framework. This probabilistic analogue is known as the Cram\'{e}r-Wold device in probability theory \cite{cramer1936some}.
\begin{theorem}[Cram\'{e}r-Wold theorem]
A distribution of a random vector $X$ in $\mathbb R^n$ is uniquely determined by the family of its push-forward distributions under the linear functionals $x \mapsto \langle v, x \rangle$, where $v \in \mathbb S^{n-1}$ is a unit vector.
\label{cramer_wold}
\end{theorem}
\begin{IEEEproof}
The first step is to relate the characteristic function of the distributions of $\langle v, X \rangle$ to that of $X$ via the simple calculation
\begin{align*}
   \varphi_{v_1 X_1 + \dots v_n X_n} (s) = \mathbb E e^{is(v_1 X_1 + \dots v_n X_n)} = \mathbb E e^{i \langle sv, X \rangle} = \varphi_X (sv).
\end{align*}
Since the left-hand side is given for all $v \in \mathbb S^{n-1}$ and all $s\in \mathbb R$, by the above identity we know the characteristic function of $X$, and thus the distribution of $X$.
\end{IEEEproof}

\begin{remark}
To see that this is in fact a probabilistic analogue of the Projection Slice Theorem, we observe that the left-hand side is simply the Fourier transform (modulo $i \mapsto -i$) of  the density of $v_1 X_1 + \dots v_n X_n$, whereas the characteristic function on the right-hand side is the Fourier transform of the joint density. The density of $v_1 X_1 + \dots v_n X_n$ on the other hand is nothing but the Radon projection orthogonal to $v$. 
\end{remark}

It is interesting to note that before the connection between tomography problems and its probabilistic counterpart was pointed out in \cite{renyi1952projections}, cf. \cite{markoe2006analytic}, the developments in both fields happened independently. The Cram\'{e}r-Wold device is used in probability theory mostly as a conceptional tool, to be more precise, it is used as a way to reduce a high-dimensional problem to a one-dimensional problem to which one can then apply well-established results. Our use of this result here is slightly different as we exploit its analogy to tomography problems explicitly to tackle the inverse problem of reconstructing the initial density from the output density.

For our inverse problem we compute the characteristic function of the output distribution to find
\begin{align}
\label{characteristic_fcn}
\begin{split}
  \varphi_{ Ce^{At} X_0 } (s) &= \mathbb E e^{\, i \,  \langle s, Ce^{At} X_0  \rangle} \\
				&= \mathbb E e^{\, i \, \langle (Ce^{At})^{\top} s, X_0 \rangle } = \varphi_{X_0} ( (Ce^{At})^{\top} s).
\end{split}
\end{align}
That is, the output distributions yield information on the characteristic function of the initial state distribution on the subspaces $\Image (Ce^{At})^{\top} = (\ker Ce^{At})^{\perp}$. This simple insight will be key in formulating characterizations for the uniqueness of the density reconstruction problem. We note however, that from \eqref{characteristic_fcn} we cannot in general gather information about the whole characteristic function due to the fact that $(Ce^{At})_{t \ge 0}$ is clearly parameterized by the scalar $t \ge 0$. This is exactly where we need to draw on analyticity properties of the characteristic function. In mathematical tomography, analyticity of the Fourier transform is typically guaranteed by the standard assumption of bounded support of the considered densities.

For the ensemble observability problem, the assumption of bounded support would exclude e.g. Gaussian distributions. We show that the assumption of bounded support can be replaced by a more general assumption for the characteristic function $\varphi_{X_0}$. We assume that for the considered initial distributions $\mathbb P_0$, the mapping $s \mapsto \varphi_{X_0}(sv) = \varphi_{\langle v, X_0 \rangle}(s)$, for all non-zero $v \in \mathbb R^n$, is real analytic, i.e. can be locally written as a power series about every point in $\mathbb R$. The role of this assumption will be further illuminated in the moment-based approach in the next section.


We begin by formulating our main result of this section, which gives a first sufficient condition for ensemble observability with respect to a specific class of initial distributions. Our result is inspired by uniqueness results for the tomographic reconstruction problem, cf. Theorem~5.2 in \cite{keinert1989inversion} and Theorem~ 3.142 in \cite{markoe2006analytic}, which give the most relaxed characterization known in the mathematical tomography literature.

\begin{theorem}
\label{thm:variety}
A linear system $(A,C)$ is ensemble observable for the class of initial distributions for which $s \mapsto \varphi_{X_0}(sv)$, for all non-zero $v \in \mathbb R^n$, is real analytic, if 
\begin{align}
 \bigcup_{t \ge 0} (\ker Ce^{At})^{\perp} =  \bigcup_{t \ge 0} \text{Im} (Ce^{At})^{\top}
\label{union}
\end{align}
is not contained in a proper algebraic subvariety of $\mathbb R^n$.
\end{theorem}

Thus, a sufficient condition is that the directions generated by $t \mapsto Ce^{At}$ are rich in the sense that we cannot find a proper algebraic variety in which the union \eqref{union} is contained in. Recall that an algebraic variety of $\mathbb R^n$ is the zero set of a polynomial, and that it is proper if it is not $\mathbb R^n$.


\begin{IEEEproof}
We show that under the analyticity condition on $\varphi_{X_0}$ and the assumption that the union~\eqref{union} is not contained in a proper algebraic variety, knowing the characteristic function on \eqref{union} is sufficient to know the characteristic function everywhere.

First of all, we consider two $\varphi_{X_0'}$ and $\varphi_{X_0''}$ such that their difference $h:=\varphi_{X_0'}-\varphi_{X_0''}$ vanishes on the union \eqref{union}, i.e.
\begin{align}
\label{eq:vanish_on_union}
 h(\xi) = 0 \; \text{ for all } \xi \in \bigcup_{t \ge 0} \text{Im} (Ce^{At})^{\top}.
\end{align}
By analyticity, we can write for any non-zero $\xi \in \mathbb R^n$ and any sufficiently small $\lambda$,
\begin{align}
h(\lambda \xi) = \sum_{p=0}^{\infty} \lambda^p a_p (\xi).
\label{eq:h_series}
\end{align}
Therein $a_p(\xi) = \frac{i^p}{p!} \left(\mathbb E(\langle \xi, X_0'\rangle^p) - \mathbb E(\langle \xi, X_0'' \rangle^p) \right)$, which is a \emph{homogeneous polynomial} of degree $p$, see Section \ref{subsec:moment_problem}.

Now for an arbitrary $\xi \in \bigcup_{t \ge 0} \text{Im} (Ce^{At})^{\top}$, by analyticity, the condition $h(\lambda \xi) = 0$, for all $\lambda$ in a neighborhood around the origin, is equivalent to the vanishing of the polynomials $a_p$ on the union \eqref{union}.
The union~\eqref{union} is thus contained in the algebraic varieties defined by $a_p$. Thus, by the assumption that the union \eqref{union} is not contained in a \emph{proper} algebraic variety, all polynomials must be trivial, i.e. $a_p \equiv 0$.

Since for all non-zero $\xi \in \mathbb R^n$, the mapping $\lambda \mapsto h(\lambda \xi)$ is real analytic in a neighborhood of any point of the real axis, $\lambda \mapsto h(\lambda \xi)$ is completely determined by its power series about the origin, which is zero. Therefore we conclude that $h \equiv 0$, i.e. $\varphi_{X_0'} = \varphi_{X_0''}$, and thus lastly $X_0'=X_0''$.
%
\end{IEEEproof}

From the main result Theorem \ref{thm:variety}, we now further derive a sufficient condition based on the special case that occurs when the affine subspaces that one is integrating over are one-dimensional. This is nowadays considered a classic result arising from the study of the X-Ray transform \cite{markoe2006analytic}. In view of the ensemble observability problem, the assumptions are quite restrictive though.

\begin{theorem}
If $(A,C)$ is observable, and $\text{rank } C= n-1$, then the union \eqref{union} is not contained in a proper algebraic variety.
\label{thm:sufficient_condition_xray}
\end{theorem}

\begin{IEEEproof}
With $\text{rank } C = n-1$, the dimension of $(\ker Ce^{At})^{\perp}$ is also $n-1$. Due to the observability of $(A,C)$, the intersection~\eqref{intersection} is trivial and thus $(\ker Ce^{At})^{\perp}$ with $t \ge 0$, constitutes an infinite family of pairwise distinct hyperplanes.  More in detail, for an observable system $(A,C)$, the mapping $t \mapsto Ce^{At}$ cannot have a discrete image $\{Ce^{A t_1}, Ce^{A t_2}, \dots \}$ for arbitrary (at most) countable times $t_1, t_2, \dots$. To see this, define the sets
\begin{align*}
     T_k := \{ t \ge 0 : Ce^{At} = Ce^{A t_k}\}.
\end{align*}
Due to continuity of $t \mapsto Ce^{At}$, the sets $T_k$ are closed. Suppose now for contradiction that $\bigcup_{k = 1, 2, \dots} T_k  = [0, \infty)$, then the Baire category theorem yields the existence of an index $k^{\star}$ so that $T_{k^{\star}}$ has non-empty interior, i.e. contains an interval. But this would contradict observability of the system during this interval.
Lastly, an infinite family of distinct hyperplanes cannot be contained in a proper algebraic variety. 
\end{IEEEproof}

Furthermore we would like to highlight the remarkable special case of $n=2$, in which the richness property is satisfied by observability of $(A,C)$ alone.
\begin{corollary}
For an observable \emph{two-dimensional} system $(A,C)$, the union \eqref{union} is not contained in a proper algebraic variety.
\label{corollary:two_dim}
\end{corollary}
%

One question that immediately comes up is whether or not an observable system $(A,C)$ already generates ``directions" rich enough such that the union \eqref{union} is not contained in an algebraic subvariety. The following example answers this question to the negative.
\begin{example}
Recall that from Theorem~\ref{thm:sufficient_condition_xray}, and Corollary \ref{corollary:two_dim}, we learned that in order to find a system that is observable, but for which the union \eqref{union} is contained in a proper algebraic variety, we need to consider systems with degree of at least three. Consider the system
\begin{align}
\begin{split}
 \dot{x}(t) &= \begin{pmatrix} 0 & 0 & 0 \\ 0 & -1 & 0 \\ 0 & 0 & -2 \end{pmatrix} x(t), \\
	y(t) &= \begin{pmatrix} 1 & 1 & 1 \end{pmatrix} x(t),
\label{variety_counter_example}
\end{split}
\end{align}
which is easily seen to be observable in the classical sense since the diagonal entries are pairwise distinct and every entry in the output matrix is non-zero. Now if we compute
\begin{align*}
  Ce^{At} = \begin{pmatrix} 1 & e^{-t} & e^{-2t} \end{pmatrix},
\end{align*}
we see that the algebraic variety given by the homogeneous polynomial equation
\begin{align}
   x_1 x_3 = x_2^2
\label{variety}
\end{align}
contains the union \eqref{union}, thus violating the richness condition in Theorem~\ref{thm:variety}.
\label{observable_but_variety}
\end{example}

Note that, since Theorem~\ref{thm:variety} is only sufficient, we can not conclude at this point that the system~\eqref{variety_counter_example} is not ensemble observable.
We will come back to this example later, after having introduced the moment-based approach to ensemble observability. There we will see that we can actually construct distinct, indistinguishable (Gaussian) initial distributions for this example.
Before that,
in the next subsection we show how the well-developed computational reconstruction methods from computed tomography can be used for the practical reconstruction of initial state distributions.

\subsection{Practical reconstruction based on tomography methods}
In this subsection, we demonstrate that the tomography framework can be also used as a means to practically reconstruct an unknown initial distribution. In particular, we also give a computational solution to the ensemble observability problem in Example \ref{example:illustrative},
\begin{align*}
 \dot{x}(t) &= \begin{pmatrix} -1 & 1 \\ 0 & 0 \end{pmatrix} x(t), \;\; x(0) \sim \mathbb P_0, \\
      y(t) &= \;\;\;\;\, \begin{pmatrix} 1 & 0 \end{pmatrix} \, x(t).
\end{align*}
Suppose that the unknown initial density is given by the bimodal density
\begin{align*}
 p_0 = 0.7 p_1 + 0.3 p_2
\end{align*}
where $p_1$ is the density of a normal distribution with mean $\mu_1 = (1,2)$ and covariance matrix $\Sigma_1 = \diag(0.3^2, 0.3^2)$ and $p_2$ is the density of a normal distribution with mean $\mu_2 = (2,1)$ and covariance matrix $\Sigma_2 = \diag(0.2^2,0.2^2)$. 

\begin{figure}[thpb]
\centering
     \includegraphics[width=0.24 \textwidth]{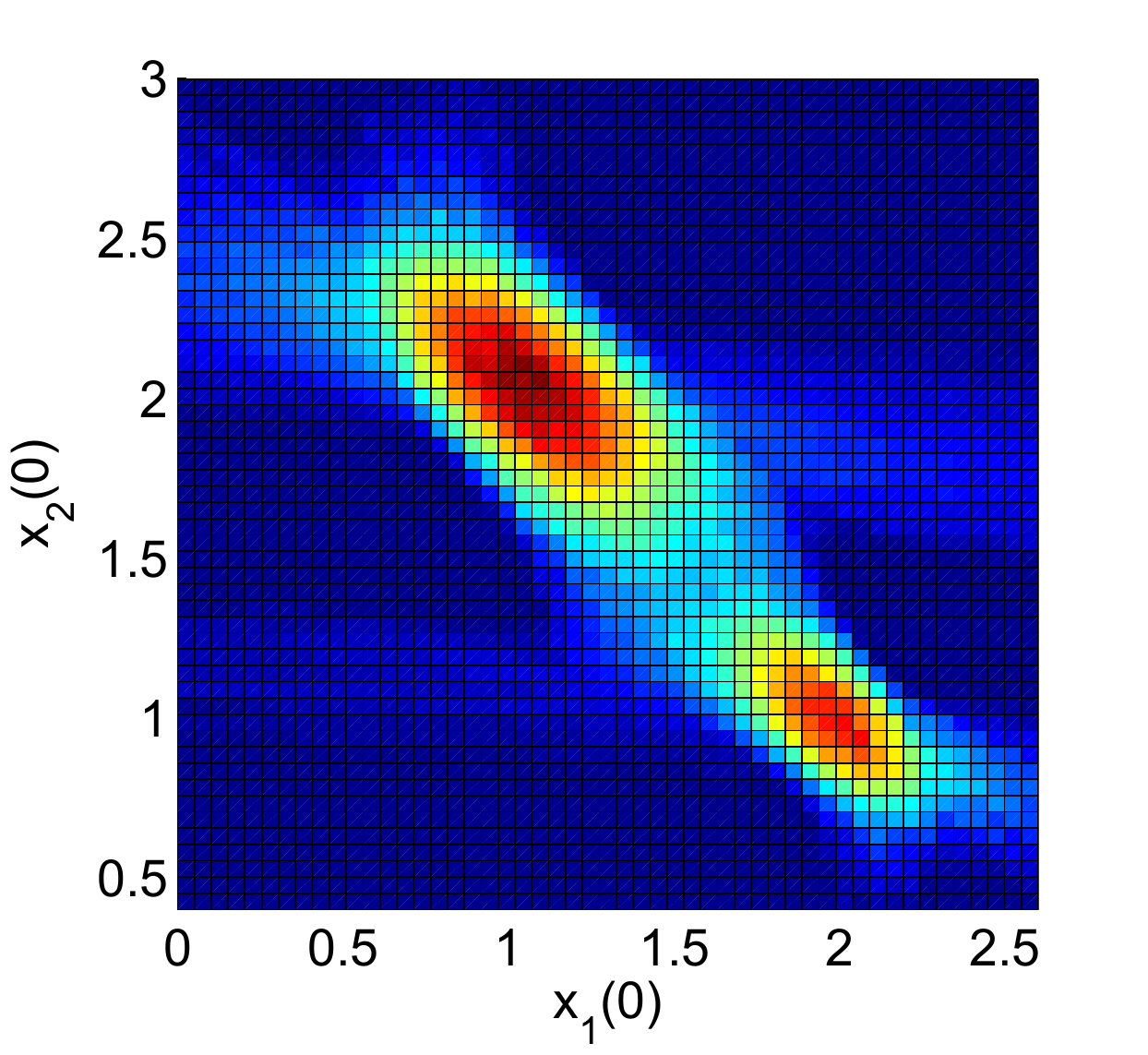}    \includegraphics[width=0.24 \textwidth]{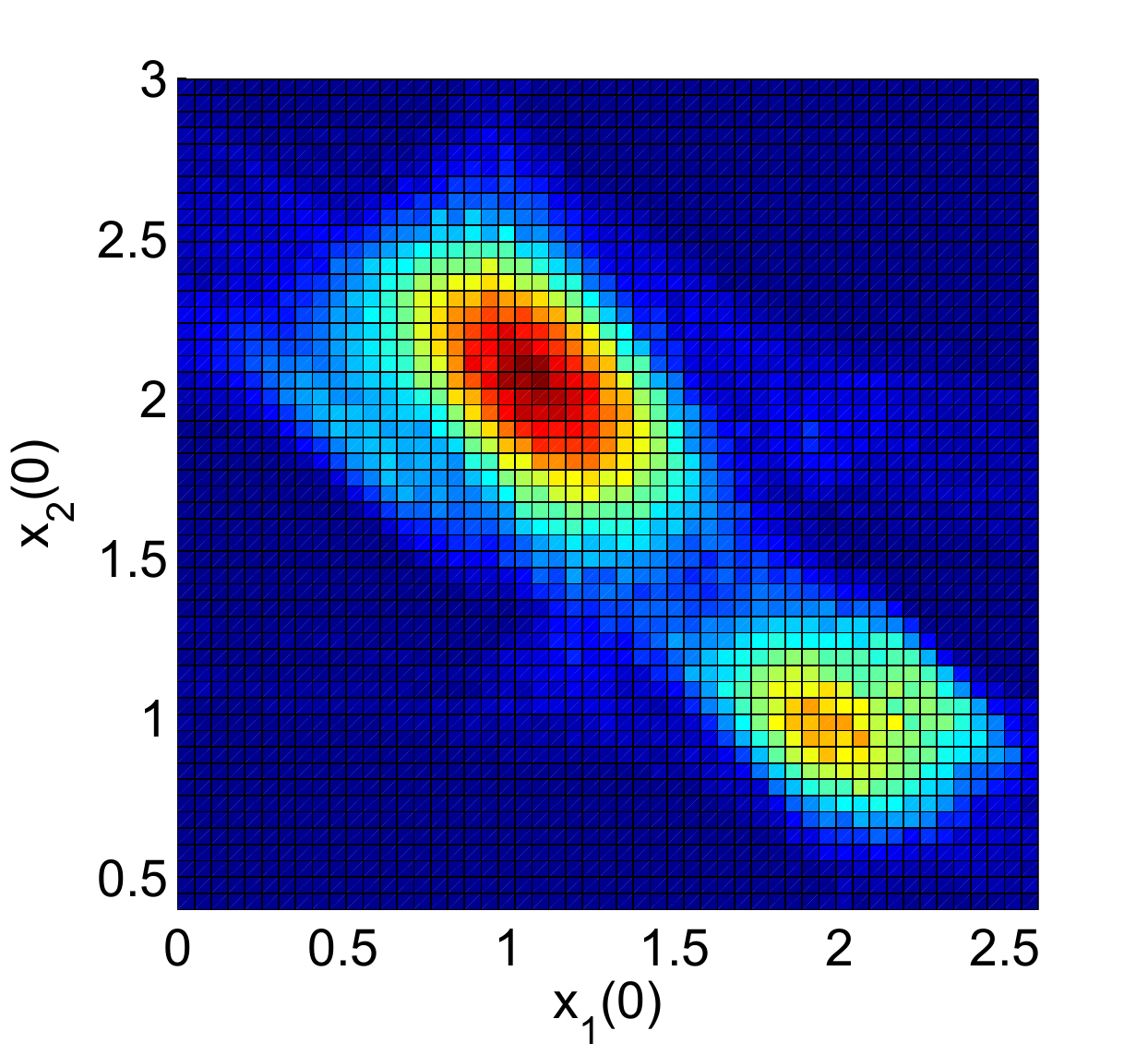}   \\ \vskip 0.2cm
      \includegraphics[width=0.24 \textwidth]{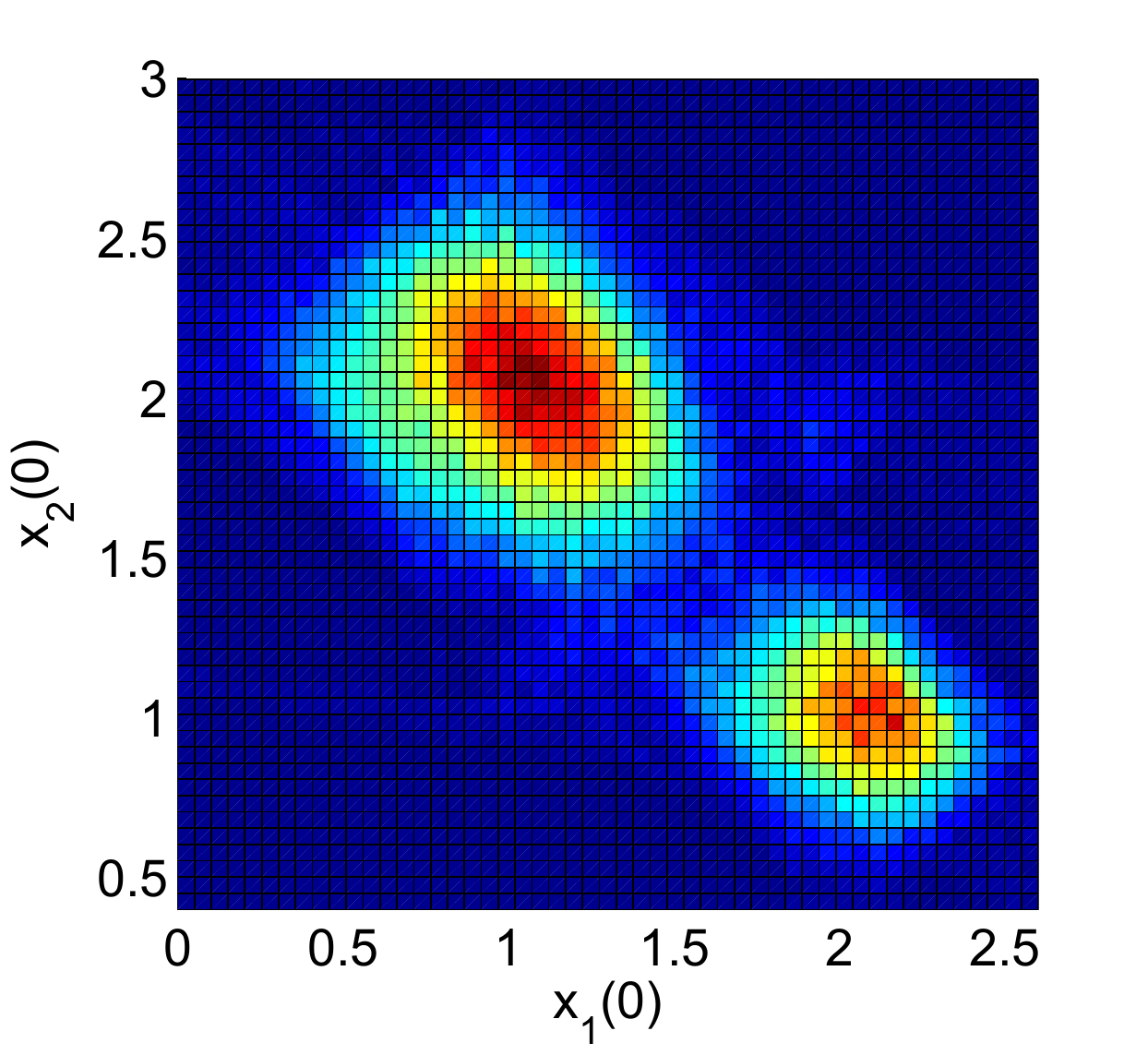}            \includegraphics[width=0.24 \textwidth]{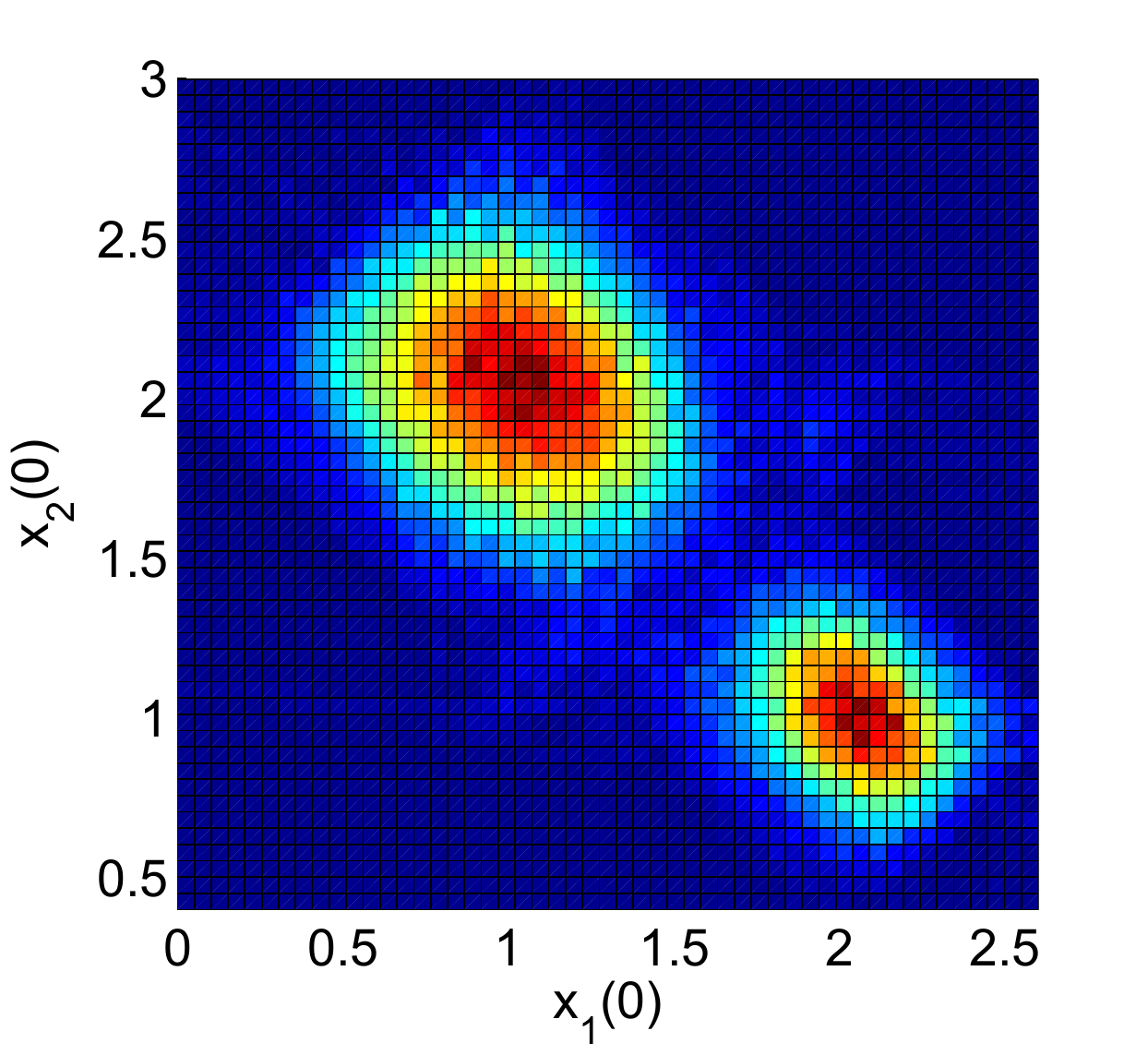}     
\caption{The reconstruction of the bimodal distribution using ART with different iteration numbers 1, 3, 5 and 7. For a smaller number of iterations, a ``distortion" effect is witnessed. This is known to be caused by the limited direction situation. By increasing the numbers of iterations, however, this effect is in large parts suppressed despite the limited direction situation.}
\label{fig:bimodal_art}
   \end{figure}

For the tomographic reconstruction of the initial density, we can employ Algebraic Reconstruction Techniques (ART), see e.g. \cite{natterer1986mathematics}. These reconstruction techniques are based on the idea of discretizing the state space into pixels so that the unknown distribution is expressed as a piecewise constant function that is constant on a fixed pixel. The strip integrals $\int_{(Ce^{At})^{-1}(B_y)} p_0 \, \text{d}x$ occuring in the tomography problem are thus approximated by weighted sums of the values of the pixels that the strip passes through. The values $\mathbb P_{y(t)}(B_y)$ are approximated via measured samples ($10^5$ samples for each time point in this example) from the output distribution. This results in a large system of linear equations which is then solved iteratively using e.g. Kaczmarz method.
The result obtained by ART is shown in Figure \ref{fig:bimodal_art} for different iteration numbers. For smaller numbers of iterations, we witness a well-known ``distortion" effect which is due to the limited direction situation, cf. Figure~\ref{fig:tilted}. The drawback of increasing the number of iterations is the numerical noise that gets introduced.

\begin{remark} While Theorem~\ref{thm:sufficient_condition_xray} combined with Theorem~\ref{thm:variety} explicitly requires measurement data of the output densities for infinitely many time points, it is clear that in practice one can only use a finite number of measured output densities. In other words, our results do not apply directly. 
To make matters worse, one can even show that there are infinitely many different densities that produce a finite set of Radon projections exactly, see e.g. \cite{markoe2006analytic}, \cite{helgason2011integral}. 
Nevertheless we can expect the practical reconstruction to yield a small estimation error for sufficiently many measured output densities.
\end{remark}

In conclusion, while the established methods of tomography seem to be generally suited for the ensemble state reconstruction problem as introduced here, there are still specific challenges stemming from the dynamic origin of this problem.
However, an exhaustive discussion of different tomographic methods for the ensemble state reconstruction problem is beyond the scope of this paper.

%

\section{Solution via observability of moments}
\label{section:moment_approach}
In this section we present an alternative, dual, route to the ensemble observability problem based on the idea of directly reconstructing the moments of the initial density. The dynamics of the moments can be described by so-called tensor systems, which are systems of homogeneous $p$-forms in the components of the state $x$ as considered by R.\ W.\ Brockett in \cite{Brockett1973_lie} and thereafter extensively studied within the observability of linear systems with polynomial output \cite{baillieul1981controllability, dayawansa1987observing, sira1988algebraic, martincarleman}. 

The problem of reconstructing the moments of the density then boils down to studying the observability of all these tensor systems which happen to be linear systems. The assumption of real analyticity of the characteristic function of the unknown density from our previous approach does in fact guarantee that moments of all orders exist and that these moments furthermore determine the density uniquely. We show as one of our main results, that both frameworks are inherently related, but that the moment-based approach in fact allows for even more general results. Having established this link, we continue with studying the ensemble observability problem in this more systems theoretic framework.

\subsection{The moment problem}
\label{subsec:moment_problem}
We start with a brief introduction to moments of multivariate probability distributions and the famous moment problem. Let $X$ be an $n$-dimensional random vector with a probability distribution $\mathbb P$. For a multi-index
\begin{align*}
 \alpha = (\alpha_1, \alpha_2, \dots, \alpha_n),
\end{align*}
i.e. an $n$-tuple of non-negative integers, we call
\begin{align*}
m_\alpha := \mathbb E[x^{\alpha}] = \mathbb E[x_1^{\alpha_1} \dots x_n^{\alpha_n}] = \int_{\mathbb R^n} x_1^{\alpha_1} \dots x_n^{\alpha_n} \, \text{d}\mathbb P(x),
\end{align*}
a moment of order $|\alpha| = \alpha_1 + \alpha_2 + \dots + \alpha_n = p$  of $\mathbb P$, if
\begin{align*}
M_p := \int_{\mathbb R^n} \| x \|^p \, \text{d}\mathbb P(x) < \infty.
\end{align*}
 A classical problem in probability theory is the question of whether or not a probability distribution is uniquely determined by its moments. This is known as the moment problem. If a distribution with finite moments is determined uniquely by its moments, the distribution is called {moment-determinate}.  

Of course, under the assumption that a probability distribution is moment-determinate, our problem of reconstructing the initial state distribution is equivalent to the reconstruction of the moments of the initial state distribution.
The close connection to our approach in the previous section is given by the observation that the mapping $s \mapsto \varphi_{X}(s v)$ is the characteristic function of the random variable $\langle v, X \rangle$. Real analyticity of $\varphi_{\langle v, X \rangle}$ on $\mathbb R$ implies that $\varphi_{\langle v, X \rangle}$ is completely determined by its power series about the origin, 
\begin{align*}
  \varphi_{\langle v, X \rangle}(s) = \sum_{p=0}^\infty s^p  \frac{ \varphi_{\langle v, X\rangle}^{(p)}(0)}{p!},
\end{align*}
cf. Section~XV.4 in \cite{feller1971introduction}. Furthermore, we have the well-known fact that the derivatives at the origin are given by
\begin{align*}
  \varphi_{\langle v, X\rangle}^{(p)}(0) = i^p \mathbb E( \langle v, X \rangle^p),
\end{align*}
i.e. $\varphi_{\langle v, X \rangle}$, due to analyticity, is uniquely determined by the moments $E( \langle v, X \rangle^p)$. These moments can be readily computed via the multinomial theorem as
\begin{align*}
 \mathbb E( \langle v, X \rangle^p) = \sum_{|\alpha| = p} \binom{p}{\alpha}  m_{\alpha} v^{\alpha},
\end{align*}
cf. the homogeneous polynomial $a_p$ in \eqref{eq:h_series}.
Now for two distributions $\mathbb P'$ and $\mathbb P''$ with the same moments, all their one-dimensional projections for non-zero $v \in \mathbb R^n$ have the same moments, and are thus equal. By virtue of the Cram\'{e}r-Wold theorem, we have $\mathbb P' = \mathbb P''$, i.e. moment-determinacy of the initial distributions.

Besides the richness condition, the approach in the previous section is essentially based on moment-determinacy of the one-dimensional projections $\langle v, X_0 \rangle$ and eventually moment-determinacy of $X_0$, which are guaranteed by the real analyticity assumption. In this section, we pursue a more direct approach which aims at directly reconstructing the moments of the initial distributions without having to go the route over analyticity arguments. In the following, we introduce the framework of so-called tensor systems in which the moment dynamics can be conveniently described.

\subsection{Background on tensor systems}
We begin by giving a brief review of tensor systems. This will also fix notation for the subsequent analysis. For a more complete introduction to tensor systems we refer to \cite{Brockett1973_lie}. Recall that for $x \in \mathbb R^n$, the vector $x^{[p]}$ denotes the vector of weighted $p$-forms in the components of $x$, i.e.
\begin{align*}
    x^{[p]} = \begin{pmatrix}  x_1^p & w_1 x_1^{p-1}x_2 & w_2 x_1^{p-1} x_3 & \dots & x_n^p \end{pmatrix}^{\top}.
\end{align*}
By a standard combinatorial ``stars and bars" argument, we conclude that the dimension of $x^{[p]}$ is $N(n,p) := \binom{n+p-1}{p}$.
More precisely, $x^{[p]}$ is the vector of weighted powers $x^{\alpha}$ with \mbox{$|\alpha|=p$}, where the entries of $x^{[p]}$ are ordered lexicographically in a decreasing order according to the multi-indices which we denote $\alpha^1, \dots, \alpha^{N(n,p)}$. In this notation, we would write
\begin{align*}
  x^{[p]} = \begin{pmatrix} w_p(\alpha^1) x^{\alpha^1} & w_p(\alpha^2) x^{\alpha^2} & \dots & w_p(\alpha^{N(n,p)}) x^{\alpha^{N(n,p)}} \end{pmatrix}^{\top},
\end{align*}
where $w_p(\alpha) := \sqrt{p!/(\alpha_1! \dots \alpha_n!)}$.

Now, if we have an equation $y=Cx$, then it can be seen that there is also a linear dependency between the vectors $y^{[p]}$ and $x^{[p]}$ which we denote by $y^{[p]} = C^{[p]} x^{[p]}$. If we consider a linear differential equation $\dot{x}(t) = A x(t)$, then it can also be seen that $x^{[p]}(t)$ satisfies a linear differential equation in which we denote the system matrix as $A_{[p]}$, i.e. $\dot{x}^{[p]}(t) = A_{[p]} x^{[p]}(t)$. Moreover, observe that the so-called tensor system
\begin{align}
\begin{split}
   \dot{x}^{[p]}(t) &= A_{[p]} x^{[p]}(t) \\
	y^{[p]}(t)     &= C^{[p]} x^{[p]}(t)
\end{split}
\label{tensor_system}	
\end{align}
has the solution $$y^{[p]}(t)  = (Ce^{At}x_0)^{[p]} = (Ce^{At})^{[p]} x_0^{[p]}.$$
Lastly, observe that by considering the transformation
\begin{align}
  x^{[p]} := \begin{pmatrix} w_p(\alpha^1) & & \\ &   \ddots & \\ & & w_p(\alpha^{N(n,p)}) \end{pmatrix} x_u^{[p]},
\label{eq:transformation_weighted_unweigted}
\end{align}
between weighted $p$-forms $x^{[p]}$ and unweighted $p$-forms $x_u^{[p]}$, observability of the LTI system for the unweighted $p$-forms is not altered. More precisely, by taking expectations in the dynamics of the unweighted $p$-forms, we obtain explicitly the dynamics of the moments of order $p$,
\begin{align*}
 \frac{d}{dt} \mathbb E[x^{[p]}] &= W^{-1} A_{[p]} W \mathbb E[x^{[p]}], \\
                              \mathbb E[y^{[p]}] &= C^{[p]} W \mathbb E[x^{[p]}],
\end{align*}
where $W$ denotes the transformation matrix in \eqref{eq:transformation_weighted_unweigted}.

\subsection{Bridging the tomography- and moment-based approaches}
In this subsection we present a bridge between the tomography-based approach with moment-based approaches based on tensor systems. We thereby establish a complete solution with two different but dual view points. The main result in this subsection is an explicit connection between the richness condition in Theorem~\ref{thm:variety}, and the unobservable subspace of a tensor system. 

\begin{lemma}
The union $\bigcup_{t \ge 0} \, \text{Im}(Ce^{At})^{\top}$ is contained in the algebraic variety defined by
\begin{align*}
   a^{\top} x^{[p]} = 0
\end{align*}
if and only if the coefficient vector $a$ is contained in the unobservable subspace of $(A_{[p]}, C^{[p]})$.
\label{thm:variety_observability1}
\end{lemma}

\begin{IEEEproof}
The condition that the union $\bigcup_{t \ge 0} \, \text{Im}(Ce^{At})^{\top}$ is contained in the algebraic variety given by $a^{\top} x^{[p]} = 0$ is equivalent to
\begin{align*}
    a^{\top}  ( (Ce^{At})^{\top} z)^{[p]} = 0
\end{align*}
for all $t \ge 0$ and $z \in \mathbb R^m$. Using $(\tilde{A} \tilde{B})^{[p]} = \tilde{A}^{[p]} \tilde{B}^{[p]}$ and 
\begin{align}
(\tilde{A}^{\top})^{[p]} = (\tilde{A}^{[p]})^{\top},
\label{eq:transpose_and_tensoring}
\end{align} 
(see e.g. \cite{Brockett1973_lie}) we arrive at
\begin{align}
   a^{\top} (Ce^{At})^{[p] \top} z^{[p]} = ((Ce^{At})^{[p]} a)^{\top} z^{[p]} =  0,
\label{eq:vanishing_polynomial_vanishing_coeff}
\end{align}
for all $t \ge 0$ and $z \in \mathbb R^m$. 

Now we recall that for a vector $v \in \mathbb R^{N(n,p)}$ the fact that
\begin{align*}
  v^{\top} z^{[p]} = 0
\end{align*}
for all $z \in \mathbb R^n$ is equivalent to $v = 0$; the vanishing of a polynomial for all values of its variables is equivalent to the vanishing of all coefficients. Thus, \eqref{eq:vanishing_polynomial_vanishing_coeff} is equivalent to $a \in \ker ((Ce^{At})^{[p]}) = \ker (C^{[p]} e^{A_{[p]}t})$ for all $t \ge 0$. In other words, the coefficient vector $a$ is contained in the unobservable subspace of the tensor system \eqref{tensor_system}, which yields the claim.
\end{IEEEproof}

\begin{remark}
The reason why weights were introduced in the beginning is exactly to guarantee that \eqref{eq:transpose_and_tensoring} is true, cf. \cite{Brockett1973_lie}. This allows us to obtain a \emph{direct} connection between coefficient vectors of the algebraic variety and unobservable states of the tensor system: Given a vector of the unobservable subspace of a tensor system, this very same vector is also the coefficient vector of a homogeneous polynomial that defines an algebraic variety in which the union \eqref{union} is contained in.
\end{remark}

Finally, we can state our main result which gives a unifying, and also more general, sufficient condition for ensemble observability based on the observability of the tensor systems.

\begin{theorem}
\label{thm:variety_observability2}
The union $\bigcup_{t \ge 0} \, \text{Im}(Ce^{At})^{\top}$ is not contained in a proper algebraic variety if and only if the systems
\begin{align*}
   \dot{x}^{[p]}(t) &= A_{[p]} x^{[p]}(t) \\
	y^{[p]}(t)     &= C^{[p]} x^{[p]}(t)	
\end{align*}
are observable for all $p \in \mathbb N$. Under these equivalent conditions, the system $(A,C)$ is ensemble observable for the class of moment-determinate initial distributions.
\label{eq:variety_coeff_unobservable_subspace}
\end{theorem}

To conclude, our main result Theorem~\ref{thm:variety_observability2} shows that the tomography approach is in fact in perfect accordance with this moment approach. 
In fact, through the systems theoretic approach we learn that in Theorem \ref{thm:variety_observability2}, moment-determinacy of the considered initial state distributions alone is sufficient, i.e. that the stronger assumption of real analyticity was in fact a technical assumption for Section \ref{section:solution_via_mathematical_tomography}. One reason for this is that the idea of the moment-based approach is more direct: Given the output distributions, we can compute their moments and then reconstruct the moments of the initial state distribution by virtue of observability of the tensor systems.

 In the following two examples, we illustrate the concepts of our theoretical framework introduced so far on the linear system considered in Example~\ref{observable_but_variety}.

\begin{example}
\label{example:covariance}
We reconsider system \eqref{variety_counter_example}. We begin with illustrating the result given in Lemma \ref{thm:variety_observability1}. First of all, we have for $x \in \mathbb R^3$ and $p=2$ the state of the second order tensor system
\begin{align*}
    x^{[2]} = \begin{pmatrix} x_1^2 & \sqrt{2} x_1 x_2 & \sqrt{2} x_1 x_3 & x_2^2 & \sqrt{2} x_2 x_3 & x_3^2 \end{pmatrix}^{\top}.
\end{align*}
Thus, for
\begin{align*}
    a = \begin{pmatrix} 0 & 0 & -\frac{1}{\sqrt{2}}& 1 & 0 & 0 \end{pmatrix}^{\top}
\end{align*}
the equation $a^{\top} x^{[2]}=0$ defines the same variety as \eqref{variety}. Furthermore it follows from $Ce^{At} = \begin{pmatrix} 1 & e^{-t} & e^{-2t} \end{pmatrix}$ and the definition of $(Ce^{At})^{[2]}$ via $(Ce^{At} x_0 )^{[2]}= (Ce^{At})^{[2]}  x_0^{[2]}$ that
\begin{align*}
 (Ce^{At})^{[2]} = \begin{pmatrix}   1 & \sqrt{2}e^{-t} & \sqrt{2}e^{-2t} & e^{-2t} & \sqrt{2} e^{-3t} & e^{-4t}    \end{pmatrix}.
\end{align*}
It is quickly verified that 
\begin{align*}
 (Ce^{At})^{[2]} a = 0,
\end{align*}
i.e. that Lemma~\ref{thm:variety_observability1} holds true. 
\end{example}

In the second example, we illuminate the meaning of the algebraic variety defined by \eqref{observable_but_variety} to the ensemble observability problem. In particular, we show that there exist observable systems $(A,C)$ which are not ensemble observable for the class of moment-determinate initial distributions.
\begin{corollary}
\label{corollary_observable_but_not_ensemble_observable}
There are systems $(A,C)$ which are observable in the classical sense, but are not ensemble observable for the class of moment-determinate initial distributions.
\end{corollary}
\begin{example}
\label{example:observable_but_not_ensemble_observable}
We consider again the observable system \eqref{variety_counter_example} and illustrate what interpretation the non-observability of the second order tensor system has in terms of moments, as well as the consequences for the ensemble observability problem. First of all, recall that we can switch between weighted and unweighted tensor vectors via the change of coordinates \eqref{eq:transformation_weighted_unweigted} given by
\begin{align*}
  x^{[p]} := \diag(1,\sqrt{2},\sqrt{2},1,\sqrt{2},1) \; x_u^{[p]}.
\end{align*}
We recall that if we take expectations in the dynamics of the unweighted monomials $x_u^{[p]}$, we get exactly the dynamics of the second order moments which is governed by the very same LTI system. Thus the transformed coordinate vector
\begin{align*}
a_u = \begin{pmatrix} 0 & 0 & -\frac{1}{2}& 1 & 0 & 0 \end{pmatrix}^{\top}
\end{align*}
is contained in the unobservable subspace of the unweighted tensor system, i.e. the system describing the second order moments. This means that if we add in the covariance matrix
\begin{align*}
\Sigma_{X_0} = \begin{pmatrix} \mathbb E[X_1^2] &  \mathbb E[X_1 X_2] &  \mathbb E[X_1 X_3] \\  \mathbb E[X_1 X_2] &  \mathbb E[X_2^2] &  \mathbb E[X_2 X_3] \\  \mathbb E[X_1 X_3] &  \mathbb E[X_2 X_3] &  \mathbb E[X_3^2] \end{pmatrix} - \mu \mu^{\top}
\end{align*}
a sufficiently small $\lambda \in \mathbb R$ to the $(2,2)$ element and $-\frac{\lambda}{2}$ to the $(1,3)$ and $(3,1)$ elements respectively such that the resulting matrix $\Sigma_{X_0}$ stays positive definite, then this will not be noticed in the output $\mathbb E[y^2 (t)]$.

In the special case that we consider Gaussian distributions, we can construct two distinct but indistinguishable initial distributions, disproving ensemble observability of system \eqref{variety_counter_example}.
%
%
For a concrete example of such indistinguishable initial state distributions, we can consider $\mathbb P_0' = \mathcal N(\mu, \Sigma')$ and $\mathbb P_0'' = \mathcal N(\mu, \Sigma'')$, with
\begin{align*}
   \Sigma' = \begin{pmatrix}   \sigma^2 & & \\ & \sigma^2 & \\ & & \sigma^2 \end{pmatrix},  \hspace{ 0.4cm }  \Sigma'' = \begin{pmatrix} \sigma^2 & 0 & - \frac{\sigma^2}{2} \\ 0 & 2 \sigma^2 & 0 \\  - \frac{\sigma^2}{2}  & 0 & \sigma^2 \end{pmatrix},
\end{align*}
i.e. we chose $\lambda = \sigma^2$. It is readily verified that for this choice, $\Sigma''$ is positive definite  (e.g. via Gershgorin circle theorem).

Now for system \eqref{variety_counter_example}, we recall that $Ce^{At} = \begin{pmatrix} 1 & e^{-t} & e^{-2t} \end{pmatrix}$, which we use for computing the covariance of the output via the well-known equation
\begin{align}
  \sigma_{y(t)}^2 = (Ce^{At}) \Sigma_{X_0} (Ce^{At})^{\top}.
\label{eq:covariance_output}
\end{align}
It is readily verified via \eqref{eq:covariance_output} that both covariance matrices $\Sigma'$ and $\Sigma''$ lead to a variance of $$(1 + e^{-2t}+ e^{-4t}) \sigma^2$$ of the output distribution. Since the output distribution is also a normal distribution, it is uniquely determined by its mean and variance, which are both the same for $\mathbb P_{y(t)}'$ and $\mathbb P_{y(t)}''$ by construction. Thus, system \eqref{variety_counter_example} which is observable in the classical sense is not ensemble observable, and in particular also gives an example for Corollary \ref{corollary_observable_but_not_ensemble_observable}.


Note that it is not immediately clear in general that the non-observability of a particular moment leads to a non-uniqueness of the reconstruction of the initial state distribution. This is because we may not be able to simply add an element of the unobservable subspace of the moment to some solution so that we have come up with a moment sequence of a different distribution that has the same output distributions. Rather, it is well-known that moments further have to satisfy additional conditions such as positive definiteness of the covariance matrix mentioned in the example. This obstructs the way for a non-uniqueness argument which is based on linearity, and shows yet again that the problem we are studying is not fully linear. 
The nonlinearity in this seemingly linear problem is naturally introduced by the consideration of probability distributions.
\end{example}

\subsection{A feasible condition for specific single-output systems}
In this subsection we derive a more practical condition on the matrices $A$ and $C$ such that for all $p \in \mathbb N$ the $p$th order tensor systems are observable. Clearly, checking the observability of infinitely many tensor systems a priori is impossible. In the following we derive a feasible necessary and sufficient condition for the observability of all tensor systems for a specific class of systems. We focus here on the class of observable single-output systems in which the system matrix has distinct real eigenvalues. The resulting condition turns out to be very restrictive, which suggests yet again that the class of ensemble observable systems must be much smaller than the class of observable systems.

\begin{theorem}
Consider an observable \emph{single-output} system, where the system matrix $A$ has distinct real eigenvalues $\lambda_1, \dots, \lambda_n$. Then the respective tensor systems are observable for all orders, if and only if the \emph{differences}
\begin{align*}
   \lambda_2 - \lambda_1, \dots, \lambda_n - \lambda_1
\end{align*}
are \emph{linearly independent over $\mathbb Q$}. 
\label{thm:linearly_independent_over_Q}
\end{theorem}

\begin{IEEEproof}
Recall that any observable single-output system with a system matrix having distinct real eigenvalues can be transformed so that
\begin{align*}
  \tilde{A} = \diag(\lambda_1, \dots, \lambda_n)
\end{align*}
with distinct real values $\lambda_i$ on the diagonal, and
\begin{align*}
  \tilde{C} = \begin{pmatrix} \tilde{c}_1 & \dots & \tilde{c}_n \end{pmatrix}
\end{align*}
where every scalar entry $\tilde{c}_i$ is non-zero, are the new system and output matrix respectively.

It is not hard to see that the matrix $\tilde{A}_{[p]}$ is then also a diagonal matrix, for which the entries on the diagonal are sums of the form $$\alpha_1 \lambda_1 + \alpha_2  \lambda_2 + \dots + \alpha_n \lambda_n,$$ where $\alpha$ is a multi-index of order $p$, i.e. $\alpha_1 + \dots + \alpha_n = p$.

Now observe that if some diagonal entries of $\tilde{A}_{[p]}$ have the same value $\tilde{\lambda}$, then in view of a Hautus test for the $p$th order tensor system, the difference $\tilde{A}_{[p]} - \tilde{\lambda}I$ will have a rank loss that is greater than one. The output matrix $\tilde{C}^{[p]}$ being a row vector however can only accomodate for exactly one rank loss.

The question of asking whether there exists $p \in \mathbb N$ so that
\begin{align*}
   \sum_{i=1}^n \alpha'_i \lambda_i =  \sum_{i=1}^n \alpha''_i \lambda_i
\end{align*}
for different multi-indices $\alpha'$ and $\alpha''$ of order $p$, can be seen to be equivalent to the question of asking whether there exists a non-zero vector $z = \begin{pmatrix} z_1 & \dots & z_n \end{pmatrix}$ of integers such that 
\begin{align*}
z_1 + \dots + z_n = 0 \text{  and  } z_1 \lambda_1 + \dots + z_n \lambda_n = 0.
\end{align*}
This is, as one can show, equivalent to $\lambda_2 - \lambda_1, \dots, \lambda_n - \lambda_1$ being linearly independent over $\mathbb Q$.
\end{IEEEproof}

The condition that real distinct eigenvalues now further have to satisfy a linear independence over the rational numbers was also reported in e.g. \cite{dayawansa1987observing,martincarleman} within the study of Carleman linearizations.

\subsection{Incorporating independence of initial state components}
In this subsection we examine what is to be gained when one has knowledge about the components of $X_0$ being independent. 
As before, our standing assumption is moment-determinacy of the considered probability distributions. 
To recap, for these distributions, the ensemble observability problem is equivalent to the reconstruction of all moments. Moreover, we assume now that the components of the random initial state $X_0$ are independent, which can be formulated as follows.

\begin{assumption}
The components $X_{0,i}$ of the random vector $X_0$ are independent, or in other words, the initial state density is decomposable as $p_0 (x) = \prod_{i=1}^n p_{0,i}(x_i)$. 
\end{assumption}

The consideration of this assumption is relevant for practical problems in which initial state distributions do satisfy such independence assumption. Furthermore, the following analysis is a nice application of our results and shows what exactly is to be gained from the independence assumption mathematically. Such independence assumption has not been considered so far in the tomography literature.

The first trick in our consideration is to consider cumulants rather than moments. It should be kept in mind that if moments exist, then cumulants do exist and that one then can directly compute moments from cumulants. Recall that for a (scalar) random variable $Y$, the cumulant-generating function is defined by
\begin{align*}
g(t) = \log \mathbb E[e^{tY}]
\end{align*}
from which the cumulants $\kappa_p$ are obtained via a power series expansion, i.e.
\begin{align*}
g(t) =\sum_{p=1}^\infty \kappa_p \frac{t^{p}}{p!}.
\end{align*}

Recall that the $p$th cumulant is homogeneous of degree $p$ in the sense that for any constant $c \in \mathbb R$,
\begin{align*}
\kappa_p(cY) = c^p \kappa_p (Y).
\end{align*}
Furthermore, if $Y'$ and $Y''$ are two independent random variables, we have additivity $\kappa_p(Y'+Y'') = \kappa_p(Y')+\kappa_p(Y'')$. The fact that independence can be naturally inserted via additivity is what makes using cumulants particularly attractive.

Now, using the independence of initial state components $X_{0,i}, \, i=1, \dots, n$, and homogeneity of the cumulants, we have for arbitrary $s \in \mathbb R^m$,
\begin{align}
\label{eq:cumulants_eq}
 \kappa_p (\langle s, Ce^{At}X_0 \rangle) =    \sum_{i=1}^n  ( (Ce^{At})^{\top} s)_i^{ \,p}  \kappa_p(X_{0,i}).
\end{align}
Therein, the left-hand side is known, and $\kappa_p(X_{0,i}), \; i=1,\dots,n$ are the unknowns that we would like to solve for. We may rewrite \eqref{eq:cumulants_eq} more compactly as
\begin{align}
 \kappa_p (\langle s, Ce^{At}X_0 \rangle) = \left\langle ((Ce^{At})^{\top}s)^{\bullet p} , \begin{pmatrix} \kappa_p (X_{0,1}) \\ \vdots \\  \kappa_p (X_{0,n}) \end{pmatrix} \right\rangle,
\label{cumulant}
\end{align}
where $\tilde{x}^{\bullet p}$ denotes the $p$th element-wise power of a vector $\tilde{x}$.

In view of our novel theoretical framework, we can see that \eqref{cumulant} can be uniquely solved for the $p$th order cumulants of $X_0$, if and only if $\bigcup_{t \ge 0} \Image (Ce^{At})^{\top}$ is not contained in a proper algebraic variety of the form
\begin{align}
   \tilde{a}_1 x_1^p + \tilde{a}_2 x_2^p + \dots + \tilde{a}_n x_n^p = 0.
\label{variety_independence}
\end{align}
Thus, by placing an independence assumption on the initial state distribution, we shrink the class of algebraic varieties to be considered in the richness condition to those algebraic varieties defined by polynomials of degree $p$, which do not have cross terms. Hereby we conveniently say that a polynomial of degree $p$ does not have a cross-term, if all monomials occurring in the polynomial are of the form
$$x^{\alpha} = x_i^{|\alpha|}.$$

Clearly this restriction makes it in principle easier for a system to be ensemble observable. From the tensor system viewpoint, $\bigcup_{t \ge 0} \Image (Ce^{At})^{\top}$ is not contained in an algebraic variety of the form \eqref{variety_independence} if and only if the \emph{intersection} of the unobservable subspace of $(A_{[p]},C^{[p]})$ \emph{with} the subspace
\begin{align}
   \big\{ a \in \mathbb R^{N(n,p)} : a_i = 0 \Leftrightarrow i\text{th entry of }x^{[p]} \text{ is a cross-term}\big\}
\label{no_mixed_constraint}
\end{align}
is trivial. This means that in view of testing observability of the tensor systems via a Hautus test, we do not need to consider every eigenvector of $A_{[p]}$, but only those that additionally lie in the set \eqref{no_mixed_constraint}. 

Although, as we pointed out, it is in principle easier to reconstruct an initial state distribution which is known to satisfy the independence assumption, it is in general \emph{also not} the case that observability of $(A,C)$ alone implies ensemble observability. We show this explicitly by constructing a concrete system in the following example.

\begin{example}
Given a three-dimensional system, it is not hard to verify that the dynamics of the second order moments
\begin{align*}
   \begin{pmatrix} \mathbb E[x_1^2] & \mathbb E[x_1 x_2] & \mathbb E[x_1 x_3] & \mathbb E[x_2^2] & \mathbb E[x_2 x_3] & \mathbb E[x_3^2] \end{pmatrix}^{\top},
\end{align*}
is described by the system matrix
\begin{align*}
\begin{pmatrix} 2 a_{11} & 2 a_{12} & 2 a_{13} & 0 & 0 & 0 \\     
				   a_{21} & a_{11}+a_{22} & a_{23} & a_{12} & a_{13} & 0 \\
				    a_{31} & a_{32} & a_{11}+a_{33} & 0 & a_{12} & a_{13} \\
				   0 & 2 a_{21} & 0 & 2 a_{22} & 2 a_{23} & 0 \\
				   0 & a_{31} & a_{21} & a_{32} & a_{22}+a_{33} & a_{23} \\
				  0 & 0 & 2 a_{31} & 0 & 2 a_{32} & 2 a_{33}  \end{pmatrix},
\end{align*}
which we denote $A_{[2]}$. Note that we dropped the normalization as it is not relevant for this analysis. The second order moment of the output $\mathbb E[y^2]$ is related to the second order moments of the state by the output matrix
\begin{align*}
  C^{[2]} := \begin{pmatrix}  c_1^2  & 2c_1 c_2 & 2c_1 c_3 & c_2^2 & 2c_2 c_3 & c_3^2 \end{pmatrix}.
\end{align*}

To recap, for the ensemble observability analysis under the independence assumption we only need to consider eigenvectors where the second, third and fifth entries are zero. To construct a counterexample, we need to find an observable $(A,C)$ failing the constrained Hautus test for the second order tensor system. In order to fail the constrained Hautus test, we need to be able to find a solution $\tilde{v}$ to the following eigenvalue problem
\begin{align*}
    A_{[2]} \tilde{v} = \begin{pmatrix}  2 a_{11} v_1 \\ a_{21} v_1 + a_{12} v_2  \\ a_{31} v_1 + a_{13} v_3 \\ 2 a_{22} v_2 \\ a_{32} v_2 + a_{23} v_3 \\ 2 a_{33} v_3  \end{pmatrix} = \lambda \begin{pmatrix} v_1 \\ 0 \\ 0 \\ v_2 \\ 0 \\ v_3 \end{pmatrix} = \lambda \tilde{v},
\end{align*}
which satisfies additionally the second condition $C^{[2]} \tilde{v} =  0$. First of all, it is seen from the eigenvalue problem that $$a_{ii} = \lambda$$ needs to hold. Now if we choose $v_1 = 1, v_2 = 1$ and $v_3=-1$, the following equations need to hold as well
\begin{align*}
   a_{21}+a_{12} = 0, \;\;\;  a_{31}-a_{13} = 0, \;\;\;  a_{32} - a_{23} = 0.
\end{align*}
Moreover, if we choose $$C = \begin{pmatrix} 1 & 1 & \sqrt{2} \end{pmatrix},$$ then $C^{[2]} \tilde{v}=0$. Based on this consideration, we come up with the system
\begin{align*}
\dot{x}(t) = \begin{pmatrix} 0 & -1 & 0 \\ 1 & 0 & 0  \\ 0 & 0 & 0  \end{pmatrix} x(t),   \hspace{0.3cm}  y(t) = \begin{pmatrix} 1 & 1 & \sqrt{2} \end{pmatrix} x(t),
\end{align*}
which can be verified to be observable, but which fails the constrained Hautus test for $(A_{[2]}, C^{[2]})$ by construction.

 For two concrete indistinguishable initial distributions, we can consider $\mathbb P_0' = \mathcal N(\mu, \Sigma')$ and $\mathbb P_0'' = \mathcal N(\mu, \Sigma'')$ with
\begin{align*}
      \Sigma' = \begin{pmatrix}   \sigma^2 & & \\ & \sigma^2 & \\ & & \sigma^2 \end{pmatrix},  \hspace{ 0.4cm }    \Sigma'' = \begin{pmatrix} \frac{3}{2} \sigma^2 & & \\ & \frac{3}{2} \sigma^2 & \\ & & \frac{1}{2} \sigma^2 \end{pmatrix},
\end{align*}
that we constructed via $\tilde{v}$ in the unobservable subspace of $(A_{[2]}, C^{[2]})$, cf. Example~\ref{example:observable_but_not_ensemble_observable}. With \eqref{eq:covariance_output} and the solution of the LTI system, $$Ce^{At} = \begin{pmatrix} \cos(t) + \sin(t) & \cos(t) - \sin(t) & \sqrt{2}\end{pmatrix},$$ we obtain for both initial distributions an output variance of $\sigma_{y(t)}^2  = 4 \sigma^2$.

Lastly, we note that Example \ref{example:observable_but_not_ensemble_observable} shows that system  \eqref{variety_counter_example} is a system which becomes ensemble observable with the additional assumption of independence. This can be seen from the unobservable subspace which is spanned by $a_u$ that is not contained in \eqref{no_mixed_constraint}. Another way to see this is through the fact that $\lambda = 0$ must hold in Example \ref{example:observable_but_not_ensemble_observable} to fullfill the independence assumption, thus obstructing us from constructing a different but indistinguishable initial distribution.

\end{example}


\section{Conclusions}

We introduced the ensemble observability problem as the problem of reconstructing the distribution of initial states in a population of finite-dimensional systems from knowledge of the output distributions over time.
The core of this problem turned out to be an inverse problem of tomographic type: the reconstruction of a density function from the knowledge of its projections which are given as integrals along a set of affine subspaces. 
This observation allowed us to apply the well-developed theory of mathematical tomography to the estimation problem.
Underlying this observation is in fact a deeper connection between observability and tomography: both problems are about inferring information about internal variables of a system, or an object, respectively, from external measurements which are projections of these variables.

The key contribution from the tomographic approach is a first sufficient condition for ensemble observability of a linear system under the assumption that the initial distributions satisfy a certain regularity property.
The tomographic approach also gives well-developed computational methods for actually solving the reconstruction problem in practical cases.
However, special characteristics of the underlying dynamic problem, such as the limited range of observation ``angles'' or the resolution of measurement times, may need more attention in numerical approaches.

While the tomographic approach gives a rather static picture of the reconstruction problem, we also pursued a more systems theoretic approach by considering the dynamics of the moments.
This led us to the study of tensor systems, as considered earlier by R.\ W.\ Brockett \cite{Brockett1973_lie}. 
Under the previously mentioned assumptions on the initial distributions, the initial distributions are fully determined by their moments.
Accordingly, ensemble observability can be inferred from the observability of all tensor systems in this case.
Furthermore, we gave a tractable condition for the observability of all tensor systems in a specific case of systems with a scalar output.


While individual systems are dynamic in this study, the ensembles are restricted to be rather static, in the sense that no individual systems are added to or removed from the ensemble.
Many ensembles in practical applications are very dynamic however; for example, cells in a cell population may divide or die.
Hence, it is also of interest to extend the concept of ensemble observability to dynamic populations.
The measure theoretic approach from which we set off in this study looks also promising for this generalization.


\bibliography{promotion_library}{}
\bibliographystyle{IEEEtran}

\end{document}